\documentclass[centertags]{amsart}
\usepackage{amsmath,amssymb,amscd}
\usepackage{hyperref}

\newtheorem*{main*}{Main Theorem}

\newtheorem{theorem}{Theorem}[section]
\newtheorem*{theorem*}{Theorem}
\newtheorem{proposition}[theorem]{Proposition}
\newtheorem{lemma}[theorem]{Lemma}
\newtheorem{corollary}[theorem]{Corollary}

\newtheorem*{question*}{Question}

\newtheorem*{conjecture*}{Conjecture}

\theoremstyle{definition}

\newtheorem{definition}[theorem]{Definition}
\newtheorem*{definition*}{Definition}

\newtheorem*{definitions*}{Definitions}

\theoremstyle{remark}
\newtheorem{remark}[theorem]{Remark}

\numberwithin{equation}{section}

\def\ds{\displaystyle}

\newcommand{\R}{\mathbb{R}}

\newcommand{\HH}{\mathbb{H}}

\newcommand{\nn}{\mathbb{N}}

\newcommand{\eps}{\epsilon }

\newcommand{\set}[1]{\left\{ #1 \right\}}
\newcommand{\op}[1]{\operatorname{#1}}

\newcommand{\dist}{\operatorname{dist}}

\newcommand{\Isom}{\operatorname{Isom}}

\newcommand{\G}{\mathsf G }

\renewcommand{\hat}{\widehat}

\newcommand{\of}{\circ }
\providecommand{\to}{\rightarrow }

\renewcommand{\bar}{\overline}

\newcommand{\abs}[1]{\left\lvert #1 \right\rvert }
\newcommand{\norm}[1]{\left\lVert #1 \right\rVert }

\newcommand{\conv}{\star}

\renewcommand{\(}{\left(}
\renewcommand{\)}{\right)}

\def\[#1\]{\begin{align*}\begin{split} #1 \end{split}\end{align*} }
\def\[[#1\]]{\begin{align}\begin{split} #1 \end{split}\end{align} }

\def\ssu{\subset}

\def\sm{\setminus}
\def\pr{\prime}

\def\nn{ \mathbb N}

\def\Ga{ \Gamma}
\def\ga{\gamma}

\def\la{ \lambda}

\def\<{\langle}
\def\>{\rangle}
\def\pr{\prime}
\def\pa{\partial}
\def\ep{\epsilon}

\begin{document}

\address{Indiana University} 
\email{connell@indiana.edu} 

\address{University of Chicago} 
\email{roma@math.uchicago.edu} 

\author[C. Connell]{Chris Connell$^\dagger$}
\thanks{$\dagger$ The author is
  supported in part by NSF grant DMS-0306594.}
\author[R. Muchnik]{Roman Muchnik$^\ddagger$}
\thanks{$\ddagger$ The author is
  supported in part by an NSF postdoctoral fellowship.}

\title{Harmonicity of Gibbs measures}

\begin{abstract}
In this paper we extend the construction of random walks with a
prescribed Poisson boundary found in \cite{ConnellMuchnik04} to the
case of measures in the class of a generalized Gibbs state. The
price for dropping the $\alpha$-quasiconformal assumptions is that
we must restrict our attention to CAT($-\kappa$) groups. Apart from
the new estimates required, we prove a new approximation scheme to
provide a positive basis for positive functions in a metric measure
space.
\end{abstract}
\maketitle


\thispagestyle{empty} 




\section{Introduction.}

In recently years a number of authors, and most notably V.
Kaimanovich, have successfully shown for various classes of
``large'' groups, that most random walks have an associated Poisson
boundary that is at least as large as its geometric boundary (e.g.
\cite{KaimanovichVershik83},\cite{Kaimanovich94,Kaimanovich00,Kaimanovich03},\cite{KaimanovichMasur96,KaimanovichMasur98}).

We are interested in the converse question of whether or not a given
measure on the geometric boundary can possibly arise as the Poisson
boundary of some random walk. Early on, Furstenberg
(\cite{Furstenberg63,Furstenberg67,Furstenberg71}) first answered
this question affirmatively for the Haar measure on the Furstenberg
boundary $G/P$ of a certain class of semisimple Lie groups $G$,
where the random walk occurs on a lattice in $G$.

In the current paper, we answer this question affirmatively for a
large class of natural measures on the geodesic boundary of a
negatively curved group. This mostly completes one direction of the
work begun in \cite{ConnellMuchnik04}, where we answered the same
question for the case of measures in the class of
$\alpha$-quasiconformal densities on the boundary of a group acting
on a Gromov hyperbolic space. The most famous example these
measures, arising for a particular $\alpha$, are the
Patterson-Sullivan classes. (In this generality, the
Patterson-Sullivan measures may possibly refer to several distinct
measure classes arising from different choices in their
construction.) We will presently treat a much larger family of
measures which include the quasiconformal ones, but in return we
must restrict the class of groups we consider.


We restrict to the case of a group $G$ acting cocompactly by
isometries on a CAT($-\kappa$) space $H$. The measure classes we
consider are represented by the {\em Gibbs streams} and arise as
certain conditional measures for the classical Gibbs states. These
classes include most of the known ergodic measures on $\pa H$. The
concept of the Gibbs state was imported into the theory of dynamical
systems directly from statistical mechanics. Gibbs streams for
negatively curved manifolds have been studied extensively in several
contexts (see \cite{Bowen75} or \cite{Ruelle78} for a list of early
references).

We briefly mention one way in which they appear in the case when $H$
is a negatively curved manifold and $G<\Isom(H)$ acts cocompactly.
Let $\mathcal{M}$ be the space of $G$ invariant Borel probability
measure on $SH$ which are invariant under the geodesic flow. Let
$h_\nu$ denote the metric entropy of the measure
$\nu\in\mathcal{M}$. If $\Phi$ is a $G$ invariant Borel function on
the unit tangent bundle $SH$, then  the {\em pressure of $\Phi$} is
given by
$$P(\Phi)=\sup_{\nu\in\mathcal{M}}\set{h_\nu-\int_{SH}\Phi
d\nu}.\footnote{The pressure of $\Phi$ was historically defined as
$P(-\Phi)$, but this choice of sign convention seems to have become
more popular and we will also follow this choice since it leads to a
more convenient normalization.}$$ A measure $\nu$ such that
$P(\Phi)=h_\nu-\nu(\Phi)$ is called an {\em equilibrium state for
$\Phi$}. If $\Phi$ is also H\"{o}lder continuous, then there is a
unique equilibrium state for $\Phi$ which is denoted by $\nu^\Phi$.
Since the geodesic flow on a closed negatively curved manifold is
Anosov, Bowen showed (\cite{Bowen75} that the equilibrium state
coincides with the Gibbs state which appears in thermodynamic
formalism as the eigenmeasure of the transfer operator corresponding
to $\Phi$. This measure has several important dynamical properties.

There is a natural identification $SH=(\pa H\times\pa H)\setminus
\Delta$ where $\Delta$ is the diagonal. We can decompose $\nu^\Phi$
into conditional measures $\set{\nu_p^\Phi}_{p\in H}$ on $\pa H$
(see \cite{Kaimanovich90}). These measures form what is called the
{\em Gibbs stream for $\Phi$}.

In Section \ref{sec:Gibbs_States}, we will exploit a construction of
U. Hamenst\"{a}dt (\cite{Hamenstadt97a}) which allows Gibbs streams
to be viewed as generalizations of $\alpha$-quasiconformal measures
in a hitherto new way. This allows us to easily generalize the Gibbs
construction to the CAT($-1$) setting.

In what follows, let $H$ be a CAT($-1$) space and let
$G<\op{Isom}(H)$ be a group of isometries acting cocompactly on $H$.
We further assume that $H$ has bounded flip (see Definition
\ref{def:flip}). This condition is satisfied for the common examples
of CAT(-1) spaces: negatively curved manifolds, buildings and
infinite trees.

Define the {\em unit tangent space} $SH$ to be the space of all unit
speed geodesics in $H$. There is a natural metric $\op{dist}$ on
$SH$ given by
$$\op{dist}(\ga_1, \ga_2) = \frac12 \int_{-\infty}^{\infty} d(\ga_1(t), \ga_2(t))
e^{-|t|} dt,$$ In this generality distinct geodesics can share
common segments, and it will be necessary to restrict our attention
to the family, $\mathcal{H}$, of those H\"{o}lder functions $\Phi$
satisfying $\Phi(\ga_1)=\Phi(\ga_2)$ if
$\ga_1([0,\eps])=\ga_2([0,\eps])$ for any $\eps>0$.


%
%
%

Now we can state our main result.
\begin{theorem}\label{thm:main1}
Let $H,G$ and $\Phi\in\mathcal{H}$ be as above. For every $p \in H$
there exists a measure $\mu_p \in P(G)$ such that $\mu_p \star
\nu_p^{\Phi} = \nu_p^{\Phi} $ and $(\partial H, \nu_p^{\Phi} )$ is
the Poisson boundary for $(G, \mu_p)$.
\end{theorem}

The above theorem is an immediate consequence of the following more
general theorem.

\begin{theorem}\label{thm:main2}
Let $H,G$ and $\Phi\in\mathcal{H}$ be as above. For every $p\in H$
and positive H\"{o}lder function $f$ on $\pa H$, there exists a
measure $\mu_p \in P(G)$ such that $\mu_p \star f\nu_p^{\Phi} =
f\nu_p^{\Phi} $ and $(\partial H, f\nu_p^{\Phi} )$ is the Poisson
boundary for $(G, \mu_p)$.
\end{theorem}

Some applications of this theorem will be explored in
\cite{BaderMuchnik05}. As an example, we indicate a consequence of a
stronger version of the above theorem which we actually prove. For
each $p,q\in G$ We may construct a mapping from $L^1(\pa G,\nu_*)$
to $l^1(G)$ given by taking $f\mapsto \set{\la_g^f(p,q)}_{g\in G}$
where $f=\sum_{g\in G} \la_g^f(p,q)\frac{d\nu^\Phi_{p
g}}{d\nu^\Phi_q}$. The important properties of this mapping are that
it is an embedding, sends positive lower semicontinuous functions to
positive sequences, and that for the standard unitary representation
$g_*f=f\of g^{-1}$ we have
$\la_\sigma^{g_*f}(p,q)=\la_\sigma^f(gp,gq)$ for all $g,\sigma\in
G$.

Note that here we have used that the right and left actions of $G$
on itself commute. So $\mu_p$ above is not $\sum_{g\in G}
\la^1_g(p,p)\delta_{gp}$ since we use $\frac{d\nu^\Phi_{p
g}}{d\nu^\Phi_p}$ in defining $\la_g^1(p,p)$ instead of the terms
$\frac{d\nu^\Phi_{g p }}{d\nu^\Phi_p}$ arising in the terms of the
convolution derivative, $\frac{d\mu_p\star
\nu^\Phi_p}{d\nu^\Phi_p}=1$.

Lastly, we would like to remark that for $G$ a $\delta$-hyperbolic
group, I. Mineyev has discovered a metric (quasi-isometric to the
original one) which mimics that of a CAT(-1) space (for instance
$d_p(x,y)=e^{-\eps (x\cdot y)_p}$ is a metric for some $\eps>0$). We
have employed a number of estimates which depend on the interior (as
opposed to asymptotic) geometry. Nevertheless, using this new metric
we suspect that it may be possible to push the main result over to
the $\delta$-hyperbolic setting. However, we did not attempt this.


\section{Regularity of measures and covers}\label{sec:regularity_covers}

We first present some necessary notation. Recall that on a set $X$,
a nonnegative function $d:X\times X\to \R$ is called a
\emph{quasimetric} (or \emph{quasidistance}) if $d$ is symmetric,
zero precisely along the diagonal, and satisfies the quasitriangle
inequality:
$$d(x,y)\leq C(d(x,z)+d(z,y))$$ for some $C\geq 1$ and all $x,y,z\in X$.

\subsection{Doubling and related properties}
 Let $X$ be a compact space equipped with a probability measure $\nu$. The symbol
$\norm{\cdot}$ will denote the $L^1(X, \nu)$-norm.

From now on, $d$ will denote either a metric or a quasi-metric on
$X$, as the context demands. Moreover, $\pi:X\times X\to [0,\infty]$
will always denote a nonnegative continuous function which is zero
precisely on the diagonal and nowhere else. We will sometimes call
such a function a distance, even though it may not be symmetric or
satisfy the triangle inequality.

Let $\Pi(x,r)=\{y\in X\,:\,\pi(y,x)< r \  {\rm and} \ \pi(x,y)<r\}$
be called an {\em open ball} with respect to such a distance, and
let $B(x,r)=\{y\in X\, :\, d(y,x)<r\}$ denote an open ball with
respect to a metric (or quasi-metric) on $X$.

\medskip

\begin{definition}\label{def:almost_decreasing}
We say that the function $f:\R\to \R_{>0}$ is {\it almost
decreasing} if there exists $\delta>0$ and $C>0$ such that
\begin{enumerate}
\item $f(x)\geq f(y)$ if $x+\delta\leq y$
\item $Cf(x) \geq f(y)$ otherwise.
\end{enumerate}
\end{definition}

\begin{definition}\label{def:nicely_decaying}
We say that an almost decreasing function
$$G: X \times X \times \R_{\geq 0}\to \R_{>0}$$ is {\it nicely
decaying with respect to the measure $\nu$} if there exists a
constant $C_G$ and $\alpha_G$, $\beta_G\geq 0$ such that
$$\sup_{x \in X}\int_{X-\Pi(x,r)}G(x,y,s)d\nu(y)\leq C_G
\frac{e^{-\alpha_G s}}{\max\{e^sr,1\}^{\beta_G}}$$
 for every $s \geq 0$ and $r \geq 0$
\end{definition}

\subsection{H\"{o}lder constants}

Recall that a map $f:X\to Y$ between metric spaces $X$ and $Y$ is
locally Lipschitz if for every $r>0$ and $x\in X$ we have,
$$\sup_{\stackrel{y \in B(x,r)}{y\neq x}}
\frac{d_Y(f(x),f(y))}{d_X(x,y)} <\infty.$$  The other extant
definitions of this notion agree when $X$ is proper. We now recall
the definition of the H\"{o}lder constant on a given scale.


\medskip

\begin{definition}\label{def:Holder_const}
As in the usual metric case, we shall say that a map $f:X\to Y$
between distance spaces $X$ and $Y$ is {\em locally
$\alpha$-H\"{o}lder} if for every $r>0$ and $x\in X$ we have,
$$\sup_{\stackrel{y \in \Pi(x,r)}{y\neq x}}
\frac{\pi_Y(f(x),f(y))}{\pi_X(x,y)^\alpha} <\infty.$$ For such a
map, we define the {\em H\"{o}lder constant at $x$ of scale $r$ and
of order $a$} to be the quantity,
$$D_r^a f(x) = \sup_{\stackrel{y \in B(x,r)}{y\neq x}}
\frac{\pi_Y(f(x),f(y))}{\pi_X(x,y)^a}.$$
\end{definition}
\medskip
 \noindent{\bf Remark: } It
is clear that if $f$ is locally $a$-H\"{o}lder and $s \leq r$ then
$D_s^af(x) \leq D_r^af(x)$. And for $b\leq a$ with $r\leq 1$ we have
$D_r^b f(x) \leq D_r^af(x)$



In the case of locally $a$-H\"{o}lder functions to $\R$, we
summarize any arithmetic relations we may need in the following
lemma. These will be mainly used in the proof of Theorem
\ref{thm:basis}. Each case may be verified by a simple (and omitted)
computation based on the definition.

\begin{lemma}\label{lem:Holder_const}
If $F, G$ are two locally $a$-H\"{o}lder functions on $X$, then
$F+G$ and $FG$ are locally $a$-H\"{o}lder. Moreover,
\begin{enumerate}
\item $D_r^a(F+G)(x) \leq D_r^a F(x) + D_r^a G(x)$,
\item $D_r^a(FG)(x) \leq \left(\sup_{d(x,y)\leq r}\abs{F(y)}\right)D_r^aG(x) + \left(\sup_{d(x,y)\leq r}\abs{G(y)}\right) D_r^a F(x)$.
\item If $G(x)\neq 0$ for all $x\in X$, then $\frac1G$ is locally
H\"{o}lder and
$$D_r\left(\frac{1}{G}\right)(x) \leq \frac{D_r^a
  G(x)}{\abs{G(x)}\left(\inf_{d(x,y)\leq r}\abs{G(y)}\right)}.$$
\item If $H:Y\to X$ is locally Lipschitz, then $F\of H$
$$D_r^a(F\circ H)(x)\leq D_{r*D_r(H)(x)}^a(F)(H(x))*[D_r(H)(x)]^a.$$
\end{enumerate}
\end{lemma}


\section{Spikes}\label{sec:spikes}
To keep the discussion as general as possible, in this section we
assume that our measure is not a single atom and has support $X$. We
will measure $X$ with respect to distance function $\pi$.

\medskip
\begin{definition}\label{def:G-spike} Assume $\G:X \times X\times \R\to \R_{>0}$
is a nicely decaying function (with respect to the measure $\nu$ and
constants $C_{\G}$ and $\alpha_{\G}$). A $5$-tuple $(h(x),r, a,s,C)$
where $h(x)$ is positive function on $X$, $r, C> 1, a \in X $, is
called a {\it ${\G}$-spike} if
\begin{enumerate}
\item $h(x)\geq \norm{h}_{L^{\infty}}/C$ on $\Pi(a, r)$, \\
\item for each $x \in B(a, r)^c$ we have
$$0<h(x)\leq C\,h(a) \,e^{\alpha_{\G} s} \int_{\Pi(a,r)} {\G}(x,y,s)d\nu(y),\quad\text{and}$$\\
\item if $y,y'\in X$ satisfy $\pi(y, y^\pr)\leq r$, then
$h(y^\pr) \leq C h(y)$.\\
\end{enumerate}

If $h(x)$ is a continuous function we call $(h(x),r, a, s, C)$ a
{\it continuous spike}. Also if $h(a)=1$ we will call $(h(x),r, a,
s,C)$ a {\it unit spike}. Lastly we will often denote the spike by
the function $h(x)$ alone with the other constants implicit.
\end{definition}
\medskip
\begin{definition}\label{def:Q-spike}  If in addition a
$\G$-spike $(h(x), r, a, s, C)$ has $h(x)$  locally $q$-H\"{o}lder
with
$$D_r^q h(x) \leq \frac{C h(x)}{r^q}$$ for all $x \in X$, 
then we call $h$ a $q$-H\"{o}lder $\G$-spike.
\end{definition}
The following two observations are immediate from the definitions.

\begin{lemma} Assume $h(x)$ is $a$-H\"{o}lder and
$$D_r^ah(x) \leq \frac{Ch(x)}{r^a},$$
then for all $r_1 \leq r$ we have $$D_{r_1}^a h(x) \leq \frac{C
h(x)}{{r_1}^a}.$$ Also a positive multiple of a spike is a spike.
\end{lemma}
%

\begin{lemma} \label{lem:comparable_spikes}Assume $(h(x), r, a, s, C)$ is a ${\G}$-spike.
If $\alpha_1 h(x) \leq f(x) \leq \alpha_2 h(x)$ for some
$\alpha_2\geq \alpha_1>0$, then $(f(x), r,a,s,
C\frac{\alpha_2}{\alpha_1})$ is a ${\G}$-spike. In particular, for
every $\alpha>0$, $(\alpha h(x), r,a, s, C)$ is still a
${\G}$-spike.\\
\end{lemma}

\begin{proof} The proof is easy and is omitted (see \cite{ConnellMuchnik04}).
\end{proof}

\section{Gibbs streams and  spikes}\label{sec:Gibbs_States}
Let $H$ be a CAT(-1)-space with metric $d$ and $\pa H$ its geodesic
boundary. The unit tangent space, $SH$, mentioned in the
introduction is usually defined to be $SH=((\pa H \times \pa H)\sm
\Delta ) \times \R$ where $\Delta$ is the diagonal subset. While
$SH$ is not a fiber bundle in general, it is homeomorphic to the
unit tangent bundle when $H$ is also a manifold. In fact, for any
fixed point $p\in H$, the homeomorphism from $SH$ to the space of
unit tangent vectors to geodesics is given by taking
$(\xi,\zeta,t)\mapsto \gamma'(0)$ where $\gamma$ is the unique unit
speed goedesic with forward endpoint $\gamma(\infty)=\xi$, backward
endpoint $\gamma(-\infty)=\zeta$ and such that $\gamma(-t)$ is the
closet point on $\gamma$ to the point $p$. (By convexity of the
metric, there is a unique closest point.) For any $CAT(-1)$ space
and choice of $p\in H$, the map carrying $(\xi,\zeta,t)\mapsto
\gamma$, with $\gamma$ defined as above, gives a bijection between
$SH$ and the space of parameterized unit speed geodesics in $H$.
Because of this, we may sometimes write $\gamma\in SH$ to indicate a
parameterized (bi-infinite) geodesic in $H$.

From now on $p,q,r$ will denote points in $H$ and $x,y, z$ will
denote points in $H \cup \pa H$. Similarly, $\zeta, \xi,$ and $\nu$
will denote points in $\pa H$.

There is a natural flow ${\mathsf g}^t$ defined on $SH$ given by
${\mathsf g}^t(\xi,\zeta,s)=(\xi,\zeta,s+t)$. Equivalently,
${\mathsf g}^t\gamma$ is the same geodesic as $\gamma$ except with
starting point shifted by $t$
 so that ${\mathsf g}^t\gamma(0)=\gamma(t)$.

We recall from the introduction, the metric $\op{dist}$ given by
$$\op{dist}(\ga_1, \ga_2)=\frac12\int_{-\infty}^{\infty} d(\ga_1(t),
\ga_2(t))e^{-|t|}dt.$$ It is relatively easy to verify that the
topology induced by $\op{dist}$ is the product topology on $((\pa H
\times \pa H)\sm \Delta ) \times \R$. By using the same formula, we
may extend the definition of $\op{dist}$ to the space $TH$ of all
(non-unit speed) geodesics in $H$, including constant geodesics
which we denote simply by their single point image.

Let $-\ga\in SH$ denote the flip of the geodesic $\ga$, i.e.,
$-\ga(t)=\ga(-t)$. We now establish some convenient properties of
$\op{dist}$.
\begin{lemma} The function $\op{dist}$ on $SH$ satisfies the following properties
\begin{enumerate}
\item $\op{dist}$ is a left invariant metric on $SH$,
\item for all $\ga\in SH$ we have
$\ds{\op{dist}(\ga,-\ga)=1},$
\item for all $\ga\in SH$ and $s\in\R$ we have $\ds{\op{dist}(\ga,{\mathsf g}^s\ga)=\abs{s}},$
\item for all $\ga_1, \ga_2\in SH$ we have
$\ds{\op{dist}(\ga_1,\ga_2)=\op{dist}(-\ga_1,-\ga_2)}.$
\end{enumerate}
\end{lemma}

\begin{proof}
The first item follows directly from the fact that $d$ is a left
invariant metric and $e^{-|t|}$ is positive. For 2) we have
$$\op{dist}(\ga, -\ga) = \frac12\int_{-\infty}^{\infty} d(\ga(t), \ga(-t))e^{-|t|}dt =\frac12\int_{-\infty}^{\infty}te^{-|t|} = 1.
$$

For 3) we have
$$\op{dist}({\mathsf g}^s \ga, \ga) = \frac12\int_{-\infty}^{\infty} d(\ga(t), \ga(t+s))e^{-|t|}dt =\frac12\int_{-\infty}^{\infty}se^{-|t|} = s.$$
Finally, we have
\begin{align*} \op{dist}(\ga_1,
\ga_2)&=\frac12\int_{-\infty}^{\infty} d(\ga_1(t),
\ga_2(t))e^{-|t|}dt =\frac12\int_{-\infty}^{\infty} d(\ga_1(-t),
\ga_2(-t))e^{-|t|}dt
\\&  =\op{dist}(-\ga_1, -\ga_2).
\end{align*}
\end{proof}

As before, let $\Phi$ be a bounded H\"{o}lder function on
$(SH,\op{dist})$. We will say $\Phi$ is {\em tempered} if whenever
$\ga_1$ and $\ga_2$ are geodesics passing through first $p$ and then
a distinct point $q\in H$,  we have $\Phi({\mathsf
g}^t\ga_1)=\Phi({\mathsf g}^t\ga_2)$ for all $0\leq t\leq d(p,q)$.
Note that such a pair of geodesics must share a common segment
between $p$ and $q$, and therefore the lift to $SH$ of a H\"{o}lder
function on $H$ is always a tempered H\"{o}lder function. However,
note that in general geodesics which intersect in a point need not
have the same $\Phi$ value.  This definition coincides with the one
given in the introduction, and we denote by $\mathcal{H}$ the space
of tempered H\"{o}lder functions.

For each distinct pair $p,q\in H$ we will denote by
$\Gamma_{p,q}\subset SH$ the set of all geodesics $\gamma_{p,q}$
such that $\ga_{p,q}(0)=p$ and $\ga_{p,q}(d(p,q))=q$. When $H$ is a
manifold $\Gamma_{p,q}$ is a single element while $\Gamma_{p,q}$ is
uncountable when $H$ is a tree, for instance. It is easy to observe
that $\ga_{p,q}((p\cdot x)_q)=\ga_{q,p}((q\cdot x)_p).$

If $\Phi$ is tempered, then we can define a new distance on $H$ by
$$d^{\Phi}(p,q) = \int_0^{d(p,q)} \Phi({\mathsf g}^t\ga_{p,q})dt.$$

\subsection{Results about $d^{\Phi}$}
\begin{lemma} Choose $\Phi\in\mathcal{H}$ with H\"{o}lder constant
$L$ and exponent $\beta$.  Then there exists a function $\hat{D}:
\R_{\geq 0} \to \R$ (depending on $L$ and $\beta$), such that for
any two geodesics $\ga_1, \ga_2$ such that $\ga_1(0)= \ga_2(0)=p$
and $d(\ga_1(T),\ga_2(T))\leq r$ we have
$$\abs{d^{\Phi}(p,\ga_1(T))-d^{\Phi}(p,\ga_2(T))} \leq \hat{D}(r).$$
\end{lemma}
\medskip
\begin{proof} Observe that
\begin{align*}
\abs{d^{\Phi}(p,\ga_1(T))-d^{\Phi}(p,\ga_2(T))} &= \abs{\int_{0}^T
\Phi({\mathsf g}^s\ga_1)-\Phi({\mathsf g}^s\ga_2)ds }  \leq
\int_0^T\abs{ \Phi({\mathsf g}^s\ga_1)-\Phi({\mathsf g}^s\ga_2)}ds
 \\ & \leq \int_0^T L\,\op{dist}( {\mathsf g}^s\ga_1, {\mathsf
g}^s \ga_2)^\beta ds \intertext{Applying Lemma
\ref{lem:integral_estimate} in Appendix \ref{sec:appendixA} we
obtain that the previous line is} &\leq L\(
\frac5{\beta}+r^{1+\beta}\).
\end{align*}

Set $\hat{D}(r)=L\( \frac5{\beta}+r^{1+\beta}\)$.
\end{proof}

\medskip
\begin{lemma} \label{lem:change_along_geod} Let $\Phi$ be a bounded function with maximum $L$. For any geodesic $\ga$ we have
$$\abs{d^{\Phi}(\ga(0), \ga(t_2))-d^{\Phi}(\ga(0), \ga(t_1))} \leq
L\abs{t_2-t_1}$$
\end{lemma}

\medskip
\begin{proof} Without loss of generality assume $t_2\geq t_1$
\begin{align*} \abs{d^{\Phi}(\ga(0),
\ga(t_2))-d^{\Phi}(\ga(0), \ga(t_1))} = \abs{\int_{t_1}^{t_2}
\Phi({\mathsf g}^s\ga)ds} \leq \int_{t_1}^{t_2} L ds =L(t_2-t_1).
\end{align*}
\end{proof}

\begin{corollary}Let $\Phi$ be a bounded H\"{o}lder function with H\"{o}lder
constant and maximum $L$.  Then for any two geodesics $\ga_1, \ga_2$
such that $\ga_1(0)= \ga_2(0)=p$ and $d(\ga_1(t_1),\ga_2(t_2))\leq
r$ for some $t_1,t_2\in \R$, we have
$$\abs{d^{\Phi}(p,\ga_1(t_1))-d^{\Phi}(p,\ga_2(t_2))} \leq \hat{D}(2r)+L\,r.$$
\end{corollary}

\smallskip
\begin{proof} Without loss of generality assume that $t_1\geq t_2$ then
we have \begin{align*}
&\abs{d^{\Phi}(p,\ga_1(t_1))-d^{\Phi}(p,\ga_2(t_2))} \leq
\abs{d^{\Phi}(p,\ga_1(t_2))-d^{\Phi}(p,\ga_2(t_2))}+\abs{d^{\Phi}(\ga_1(t_2),
\ga_1(t_1))}. \intertext{By convexity of the metric,
$\abs{t_1-t_2}\leq r$ and therefore
$$d(\ga_1(t_2), \ga_2(t_2)) \leq d(\ga_1(t_1), \ga_2(t_2))
+d(\ga_1(t_1), \ga_1(t_2)) \leq 2r.$$ And so we have}
&\abs{d^{\Phi}(p,\ga_1(t_1))-d^{\Phi}(p,\ga_2(t_2))}  \leq
\hat{D}(2r) + L\,r.
\end{align*}
\end{proof}

\begin{corollary}\label{cor:compare_geodesics} Let $\Phi$ and $\hat{\Phi}$ be H\"{o}lder functions with H\"{o}lder
constants and maximum $L$. There exists a function $D: \R_{\geq 0}
\to \R$, such that for any two geodesics $\ga_1$ and $\ga_2$ with
$d(\ga_1(0), \ga_2(0)) \leq r$ and $d(\ga_1(t_1), \ga_2(t_2)) \leq
r$ for some $t_1,t_2\in \R$ we have
$$\abs{d^{\Phi}(\ga_1(0), \ga_1(t_1))-d^{\Phi}(\ga_2(0),\ga_2(t_2))}
\leq D(r).$$
\end{corollary}

\medskip
\begin{proof} Let $\ga_3$ be a geodesic such that $\ga_3(0)=\ga_1(0)$
and $\ga_3(t_3)=\ga_2(t_2)$.  Then \begin{align*}
&\abs{d^{\Phi}(\ga_1(0), \ga_1(t_1))-d^{\Phi}(\ga_2(0),\ga_2(t_2))}
\leq
\\& \leq \abs{d^{\Phi}(\ga_1(0),
\ga_1(t_1))-d^{\Phi}(\ga_3(0),\ga_3(t_3))}+ \abs{d^{\Phi}(\ga_3(0),
\ga_3(t_3))+d^{\Phi}(\ga_2(0),\ga_2(t_2))} \intertext{Since
$d^{\Phi}(p,q)= d^{\hat{\Phi}}(q,p)$, $\ga_3(0)=\ga_1(0)$ and
$\ga_3(t_3)=\ga_2(t_2)$, we have} &\abs{d^{\Phi}(\ga_1(0),
\ga_1(t_1))-d^{\Phi}(\ga_2(0),\ga_2(t_2))}\leq 2\left(\hat{D}(2r) +
Lr\right).
\end{align*}
Set $D(r)=2\hat{D}(2r)+2L\,r$. \end{proof}

\subsection{$\Phi$-Busemann function}

For any  H\"{o}lder function $\Phi$ we define its corresponding
$\Phi$-Busemann function as follows. For fixed $p,q \in H$, we have
$$\rho^{\Phi}_x(p,q) = \lim_{z \to x} (d^{\Phi}(q,z) - d^{\Phi}(p,z)).$$
If $\Phi$ is symmetric we have $\rho^{\Phi}_x(p,q) =
-\rho^{\Phi}_x(q,p)$. When $\Phi\equiv 1$ then $d^\Phi$ is the
original metric $d$ and hence $\rho^\Phi_x=\rho^1_x$ is just the
ordinary Busemann function, which we may sometimes denote by
$\rho_x$.

Recall that the Gromov-product is defined as
$$(x\cdot y)_p = \frac{1}{2}(d(x,p)+d(y,p)-d(x,y)),$$
where in the case when either $x\in \pa H$ or $y \in \pa H$, we take
the corresponding limits (see \cite{Coornaert93}). It is easy to
observe that
$$2(q\cdot \xi)_p=\rho_{\xi}(q,p)+d(p,q).$$

Recall that any $CAT(-1)$-space is $\delta$-hyperbolic for a
sufficiently small value of $\delta$.
\begin{lemma}\label{lem:approx_busemann} Choose $\delta$ so that $H$ is $\delta$-hyperbolic. For any bounded H\"{o}lder function $\Phi$,
we have
$$\abs{\rho^{\Phi}_x(p,q)+ d^{\Phi}(p,\ga_{q,p}((q\cdot x)_p)) -
d^{\Phi}(q,\ga_{p,q}((p\cdot x)_q))}\leq 3D(\delta).$$
\end{lemma}
\begin{proof} Recall that for any triangle $p,q,x$ in
$\delta$-hyperbolic space we have
$$\op{diam}\set{\ga_{p,q}((q\cdot x)_p), \ga_{p,x}((q\cdot x)_p),
\ga_{q,x}((p\cdot x)_q)} \leq \delta.$$ So
\begin{align*}
&\abs{\rho^{\Phi}_x(p,q)+ d^{\Phi}(p,\ga_{q,p}((q\cdot x)_p)) -
d^{\Phi}(q,\ga_{p,q}((p\cdot x)_q))}\\& = \abs{\lim_{z \to
x}(d^{\Phi}(q,z)-d^{\Phi}(p,z))+d^{\Phi}(p,\ga_{q,p}((q\cdot x)_p))
- d^{\Phi}(q,\ga_{p,q}((p\cdot x)_q))} \\ &=
 \Biggr|\lim_{z \to
x}d^{\Phi}(\ga_{q,x}((p\cdot x)_p),z)-d^{\Phi}(\ga_{p,x}((q\cdot
x)_p),z)\\ &\hspace{2cm}+(d^{\Phi}(p,\ga_{q,p}((q\cdot
x)_p))-d^{\Phi}(p, \ga_{p,x}((q\cdot x)_p)))\\ &\hspace{2cm}+
(d^{\Phi}(q,\ga_{q,x}((p\cdot x)_p))-d^{\Phi}(q,\ga_{p,q}((p\cdot
x)_q)))\Biggr|  \\ & \leq 3 D(\delta).
\end{align*}
The last inequality follows from Corollary
\ref{cor:compare_geodesics}.
\end{proof}

\subsection{Patterson-Sullivan
construction}\label{sub:Patterson-construct}

Let $G$ be a subgroup of $\Isom(H)$. Recall that the {\it Poincare
series} associated to $G$ is
$$S^{\Phi}(\la)= \sum_{g \in G} e^{-d^{\Phi}(p, gp)- \la d(p, gp)},$$
for a choice of fixed point $p \in H$. Define the {\em critical
exponent} for $\Phi$ as
$$\la_{\Phi}=\sup_{\la}\{\la: S^{\Phi}(\la)=\infty\} = \inf_{\la}\{\la:S^{\Phi}(\la)<\infty\}.$$
This is independent of the choice of $p$.  It is easy to observe
that for any constant $C$ we have $\la_{\Phi +C}= \la_{\Phi}-C$ and
$\mu_p^{\Phi+C} =\mu_p^{\Phi}$.

In the event that $S^\Phi(\la_\Phi)<\infty$, we can repeat the above
construction with
$$S^{\Phi}(\la)= \sum_{g \in G} h(d^{\Phi}(p, gp)+\la d(p,gp))e^{-d^{\Phi}(p, gp)- \la d(p, gp)},$$
where $h:\R^+\to\R^+$ is a suitable subexponentially increasing
function. Patterson \cite{Patterson76} showed that $h$ can be chosen
so that the critical exponent for the modified Dirichlet series
remains $\la_\Phi$ and $S^\Phi$ diverges at $\la=\la_\Phi$. Note
that this is a general statement about series of real numbers of the
form $\sum_{i=1}^\infty h(x_i)e^{-x_i}$ for a positive sequence
$\set{x_i}$.

For any $q\in H$, the measure $\mu_q^{\Phi}$ is the $weak^{\star}$
limit of the normalized measures
$$\frac{1}{S^{\Phi}(\la)}\sum_{g \in G} h(d^{\Phi}(p, gp)+\la d(p,gp)) e^{-d^{\Phi}(q, gq)- \la
d(q, gq)} \delta_{gq},$$ as $\la \searrow \la_{\Phi}$. Since $G$
acts cocompactly on $H$, the measure $\mu_q^\Phi$ is independent of
the choice of limiting sequence in $\lambda$. The proof is the same
as that for the uniqueness of conformal densities in CAT(-1) spaces
(see \cite{Coornaert95}). By construction $\mu^\Phi_q$ is a finite
measure supported on $\pa H$. It is a probability measure when
$q=p$, but in general it is only finite for other points $q$.

\medskip
From now on $\hat{\Phi}$ denotes the flip of $\Phi$. Also define
$Sym{(\Phi)}=\frac{\Phi +\hat{\Phi}}{2}$.
\begin{lemma}
If $\Phi$ is $G$-invariant, then $$\la_{\Phi}=\la_{\hat{\Phi}} \geq
\la_{Sym{(\Phi)}}.$$
\end{lemma}

\begin{proof} It is easy to observe that \begin{align*}d^{\Phi}(p,g p)=
&\int_0^{d(p,g p)} \Phi({\mathsf g}^t\ga_{p,gp})dt =\\ &=
\int_0^{d(p,g p)} \Phi({\mathsf g}^t\ga_{g^{-1} p,p})dt
=\int_0^{d(p,g p)} \hat{\Phi}({\mathsf g}^t\ga_{p,g^{-1}p})dt
=d^{\hat{\Phi}}(p, g^{-1} p).
\end{align*}

This implies that $\la_{\Phi}=\la_{\hat{\Phi}}$. Also applying the
inequality $a^2+b^2\geq 2ab$ to $a=e^{-d^\Phi(p,gp)}$ and
$b=e^{-d^\Phi(p,g^{-1}p)}$, we obtain
$$S^{\Phi}(\la)+S^{\hat{\Phi}}(\la) \geq
2S^{Sym(\Phi)}(\la).$$

 So we obtain $\la_{\Phi}=\la_{\hat{\Phi}} \geq
 \la_{Sym(\Phi)}.$ \end{proof}

Define $\Delta(\Phi)= \la_{Sym(\Phi)} -\la_{\Phi}\leq 0$, so that
$\la_{\Phi}=\la_{Sym(\Phi)+\Delta(\Phi)}$.

The next proposition appears in \cite{Kaimanovich04} as Proposition
3.5 for the case of manifolds. However, the measures, $\nu_p^\Phi$,
that the author defines relate to the ones we define, $\mu_p^\Phi$,
by $\nu^\Phi_p=\mu^{-\Phi}_p$. We choose this normalization for
convenience later on. For the general CAT(-1) case, the proof is a
simple calculation from the above formula for $\mu_p$ and is
identical to the original case found in \cite{Sullivan79}.

\begin{proposition} For $G<\Isom(H)$ and any $G$-invariant
H\"{o}lder function $\Phi$ on $SH$, there exists a family of finite
positive measures $\{\mu_p^{\Phi}\}_{p \in \pa H}$ parametrized by
points $ p \in H$ with the following two properties:
\begin{enumerate}
\item The family $\{\mu_p^{\Phi}\}$ is $G$-equivariant, i.e.,
$\mu_{gp}^{\Phi}=g \mu_p^{\Phi}$ for any $p \in H$ and $g \in G$;

\item  The measures $\mu_p^{\Phi}$ are pairwise equivalent and
$$\frac{d\mu_q^{\Phi}}{d\mu_p^{\Phi}}(\xi)=e^{-\rho^{\Phi}_{\xi}(p,q)-\la_\Phi \rho_\xi(p,q)}.$$
\end{enumerate}
The family $\{\mu_p^{\Phi}\}$ with these properties is unique up to
global constant multiple.
\end{proposition}

\begin{remark}
As remarked in the introduction, the critical exponent $\la_\Phi$
coincides with the topological pressure of $\Phi$ (see Kaimanovich
\cite{Kaimanovich04}, but note that he is using $P(-\Phi)$ for
$P(\Phi)$). For instance, when $\Phi=C$, the formula for $\la_\Phi$
given earlier, above implies that its pressure is
$\la_{C}=\lambda_0-C$ where $\la_0$ is the ordinary critical
exponent of the group $G$ relative to the metric $d$. On the other
hand, the variational formula for the pressure of $0$, implies that
$\la_0$ is the topological entropy of the geodesic flow. At least in
the case when $H$ is a manifold and $G$ acts cocompactly, the
pressure functional defined in the introduction also agrees with the
Bowen-Ruelle definition arising from thermodynamic formalism (see
Ruelle \cite{Ruelle78}).
\end{remark}

\subsection{The Function ${\mathsf G}$}

We now continue our analysis of $d^\Phi$.

\begin{definition}\label{def:geodesic_average}
For $\eps\in \R$ and a point $p\in H$, we say that a function $\Phi$
has {\it geodesic average at least $\ep$ at $p$} if there is a
constant $T\in\R$ such that
$$d^{\Phi}(p, \ga_{p,\xi}(s)) \geq s\,\ep-T,$$
for all $s\geq 0$ and $\xi\in \pa H$.
\end{definition}

Note that the notion of geodesic average is stronger than the
ordinary time average along a geodesic since we assume $T$ is a
fixed value instead of a sublinear function in $s$.

\medskip

Now we define ${\mathsf G}^{\Phi}:\pa H \times \pa H \times \R_{\ge
0} \to \R_{>0}$ as follows
\begin{gather}
{\mathsf G}^{\Phi}(x,y, s) = \begin{cases}
e^{-2d^{\Phi}(\ga_{p,y}((x\cdot y)_p),\ga_{p,y}(s))} & \text{ if } s\geq(x\cdot y)_p \\
1 & \text{otherwise}
\end{cases}
\end{gather}

Let $\mu_p$ be Patterson-Sullivan measure at the point $p$, i.e.,
the Gibbs stream based at $p$ for the function $\Phi=0$ with
corresponding critical constant $\la_0$.
\begin{lemma}   Assume that $\Phi-\la_{\Phi}$ has positive
geodesic average $\ep$ at $p$ for $\ep>0$.  Then
$$\int_{\{y\,:\,(x\cdot y)\geq a\}} {\mathsf
G}^{\Phi}(x,y,s)d\mu_p(y) \leq C_{\mu_p}e^{-\la_0
s}\min\{e^{-\ep(s-a)},1\} .$$
\end{lemma}

\begin{proof}
By compactness, we may assume  $\abs{\Phi}\leq L$ where $L$ is the
H\"{o}lder constant for $\Phi$. Let $0=t_0<t_1<\dots <
t_{n+1}=\min{(a,s)}$ such that $t_{i+1}-t_i=1$. Let $\Pi_i=\{y \in
\pa H: t_i\leq(x\cdot y)_p\leq t_{i+1}\}$. By Lemma
\ref{lem:change_along_geod} for all $s, t$ such that $\abs{s-t}\leq
1$ we have $\abs{2d^{Sym(\Phi)}(\ga_{p,x}(s), \ga_{p,x}(t))}\leq
2L.$
\begin{align*}
\int_{\cup \Pi_i} {\mathsf G}^{\Phi}(x,y,s)d\mu_p(y) &\leq
\sum_{i=0}^{n} \int_{\Pi_i} e^{-2d^{\Phi}(\ga_{p,y}((x\cdot
y)_p),\ga_{p,y}(s))}d\mu_p(y)
\\& \leq \sum_{i=0}^{n} \int_{\Pi_i}
e^{2L}e^{-2d^{\Phi}( \ga_{p,y}(t_i),\ga_{p,y}(s))}d\mu_p(y)
\intertext{ Since $2 d^{f}=d^{f-g}+d^{f+g}$ we obtain} &\leq
\sum_{i=0}^{n} \int_{\Pi_i} e^{2L}e^{-d^{\Phi-\la_{\Phi}}(
\ga_{p,y}(t_i),\ga_{p,y}(s))}e^{-d^{\Phi+\la_{\Phi}}(
\ga_{p,y}(t_i),\ga_{p,y}(s))}d\mu_p(y)
\\
&\leq \sum_{i=0}^{n} e^{2L}e^{-\ep (s-t_i)}\int_{\Pi_i}
e^{-d^{\Phi+\la_{\Phi}}( \ga_{p,y} (t_i),\ga_{p,y}(s))}d\mu_p(y)
\intertext{since $\rho_{\xi}^0(p,\ga_{p,x}(t_i)) \leq -t_i +1+
2\delta$, for all $\xi$ such that $(\xi\cdot x)_p \geq t_i$ where
$\pa H$ is $\delta$-hyperbolic space.  Hence
$\frac{d\mu_{\ga_{p,x}(t_i)}}{d\mu_p}(\xi)=
e^{-\la_0\rho_{\xi}^0(p,\ga_{p,y}(t_i))}\geq
e^{\la_0(t_i-1-2\delta)}$, yielding} & \leq
\sum_{i=0}^{n}e^{2L}e^{-\ep (s-t_i)} \int_{\Pi_i}
e^{-d^{\Phi+\la_{\Phi}}( \ga_{p,y}(t_i),\ga_{p,y}(s))}
e^{-\la_0(t_i-1-2\delta)}d\mu_{\ga_{p,x}(t_i)}(y)\\
\intertext{since $d(\ga_{p,x}(t_i), \ga_{p,y}(t_i))\leq \delta+2$,
because $t_i+1\geq (x\cdot y)_p \geq t_i$.  So we have
$d(\ga_{\ga_{p,x}(t_i),y}(s-t_i),\ga_{p,y}(s))\leq 2+2\delta$. From
this estimate we obtain the bound, $\abs{d^{\Phi}(\ga_{p,x}
(t_i),\ga_{p,y}(s)) - d^{\Phi}(\ga_{p,x}(t_i),
\ga_{\ga_{p,x}(t_i),y}(s-t_i))}\leq D(2+2\delta)$. Now by Corollary
\ref{cor:compare_geodesics} we can continue from the previous
inequality,}
&\leq \sum_{i=0}^{n}e^{2L}e^{-\la_0(t_i+1+2\delta)} e^{-\ep (s-t_i)}
\\ & \phantom{blah blah }
\cdot\int_{\pa H} e^{-d^{\Phi+\la_{\Phi}}(\ga_{p,x}(t_i),
\ga_{\ga_{p,x}(t_i),y}(s-t_i))+D(2+2\delta)}
d\mu_{\ga_{p,x}(t_i)}(y)
\\& \leq \sum_{i=0}^{n}e^{2L+D(2+2\delta)}e^{-\la_0(t_i+1+2\delta)} e^{-\ep (s-t_i)}
e^{-\la_0 (s-t_i)}
\\ &\leq C^{\pr} e^{-\la_0 s}e^{-\ep(s-\min{(s,a)})}
\end{align*}
To complete the above estimate we need to show that for all $q \in
H$, for all H\"{o}lder functions $\Phi$  and $s\geq 0$, we have
$$\int_{\pa H} e^{-d^{\Phi+\la_{\Phi}}(q,\ga_{q,\xi}(s))} d\mu_q(\xi) \leq
e^{-\la_0 s}.$$ This is precisely the statement of Corollary
\ref{cor:shadow} in the Appendix C, since $\la_{\Phi+\la_\Phi}=0$.
If $a\leq s$ then the proof is done. In case of $a \geq s$ observe
that $\{y\in \pa H\,:\,s \leq (x\cdot y)_p \leq a\} \ssu \{y\in \pa
H\,:\,s \leq (x\cdot y)_p <\infty\} =\Pi_{\infty}$ and
$$\int_{\Pi_{\infty}} {\mathsf G}^{\Phi}(x,y,s)d\mu_p(y) \leq \mu_p(\Pi_{\infty}) \leq
C^{\pr\pr}e^{-\la_0 s}.$$
\end{proof}

We now recall the definition of ``nicely decaying" from
\ref{def:nicely_decaying}.
\begin{corollary} If $\Phi - \la_{\Phi}$ has a positive geodesic
average $\ep$ at $p$, then the function ${\mathsf G}^{\Phi}$ is
nicely decaying with respect to the Patterson-Sullivan measure
$\mu_p$ with constants $\alpha_G=\la_0$ and $\beta_G=\ep$.
\end{corollary}

Since $\Phi+C-\la_{\Phi+C}=\Phi-\la_{\Phi}+2C$, we may always
arrange, without changing $\nu_p^\Phi$, that $\Phi - \la_{\Phi}$
have positive geodesic average by simply adding a sufficiently large
constant to $\Phi$. Moreover, in the case when $G$ acts cocompactly,
Proposition \ref{prop:positive_average} in Appendix
\ref{app:decay_cond} shows that when $\Phi$ is normalized to have
zero pressure it will automatically have a uniformly positive
average along all geodesics.

\subsection{Derivatives are ${\mathsf G}^{Sym(\Phi)}$-spikes} For the purpose of this
section we take $L$ to be a number larger than the H\"{o}lder
constants for $\Phi$ and $\hat{\Phi}$ and the upper bound for
$\abs{\Phi}$. As a reminder in what follows, recall that $\ga_{p,q}$
represents {\em any} ray which extends the oriented geodesic segment
from $p$ to $q$. All results are independent of this choice, when
there is a choice.

\begin{lemma} Assume  $\Phi\in\mathcal{H}$  has positive geodesic
average (with constants $\ep>0$ and $T$) at some $p \in H$.  There
exists a constant $C_{\Phi}$ such that for all $q \in H$, the
$5$-tuple $$\(e^{-2d^{\Phi}( \ga_{p, q}((\xi\cdot q)_p),q)},
\ga_{p,q}(\infty), e^{-d(p,q)}, d(p,q),C_{\Phi}\)$$ is a ${\mathsf
G}^{\Phi}$-spike with respect to the symmetric distance
$\pi_p(x,y)=e^{-(x\cdot y)_p}$ (actually a metric after possibly
rescaling $d$).
\end{lemma}

\medskip
\begin{proof}  Set $h(\xi)=e^{-2d^{\Phi}( \ga_{p, q}((\xi\cdot
q)_p),q)}$. It is clear that $h(\ga_{p,q}(\infty)) =1 \geq
\max(h)e^{-2T}$.

Observe that, for all $\xi, \zeta \in \pa H$ such that
$(\xi\cdot\zeta)_p\geq d(p,q)$, we have $\abs{(\xi\cdot q)_p -
(\zeta\cdot q)_p}\leq \delta$.  Thus
$$\frac{h(\xi)}{h(\zeta)} = e^{2d^{\Phi}(
\ga_{p, q}((\zeta\cdot q)_p),q)-2d^{\Phi}( \ga_{p, q}((\xi\cdot
q)_p),q)}\leq e^{\abs{2d^{\Phi}( \ga_{p, q}((\zeta\cdot
q)_p),q)-2d^{\Phi}( \ga_{p, q}((\xi\cdot q)_p),q)}} \leq
e^{2L\delta}.$$

Recall that $\Pi(\ga_{p,q}(\infty),e^{-d(p,q)}) = \{y \in \pa H\,:
(y\cdot \ga_{p,q}(\infty))_p \geq d(p,q)\}$. Also
$\mu_p(\Pi(\ga_{p,q}(\infty), e^{- d(p,q)}) \geq
\frac{1}{K}e^{-\la^0 d(p,q)}$. So we have \begin{align*}h(x) &=
\frac{1}{\mu_p(\Pi(\ga_{p,q}(\infty),e^{-
d(p,q)}))}\int_{\Pi(\ga_{p,q}(\infty),e^{- d(p,q)})} h(x) d\mu_p(y)
\\ & \leq h(\ga_{p,q}(\infty)) Ke^{\la^0
d(p,q)}\int_{\Pi(\ga_{p,q}(\infty),e^{-d(p,q)})} h(x) d\mu_p(y).
\end{align*}

To finish the proof, we need to show that
$\frac{h(x)}{G(x,y,d(p,q))}$ is uniformly bounded for every $y \in
\Pi(\ga_{p,q}(\infty),e^{- d(p,q)})$. We need to consider $2$ cases:

Case 1) (Easy Case) $d(p,q)\leq (x\cdot y)_p$. Then $h(x) \leq
e^{2T} \leq e^{2T}G(x,y,d(p,q))$.

 Case 2) $d(p,q) \geq (x\cdot y)_p$.  In this case
$G(x,y,d(p,q))= e^{-2d^{\Phi}(\ga_{p,y}((x\cdot
y)_p),\ga_{p,y}(d(p,q)))}$.

Since $(y\cdot \ga_{p,q}(\infty))_p \geq d(p,q)$, we have $(y\cdot
q)_p \geq d(p,q)-\delta.$  So
$$(x\cdot y)_p \geq \min((x\cdot q)_p, (y\cdot
q)_p)-\delta \geq \min((x\cdot q)_p, d(p,q)-\delta)-\delta \geq
(x\cdot q)_p-2\delta.$$ On the other hand, $$(x\cdot q)_p \geq
\min((x\cdot y)_p, (y\cdot q)_p ) - \delta  \geq \min((x\cdot y)_p,
d(p,q)-\delta)- \delta \geq (x\cdot y)_p - 2\delta.$$ We obtain
$\abs{(x\cdot y)_p - (x\cdot q)_p}\leq 2\delta$.

Since $(y\cdot \ga_{p,q}(\infty))_p \geq d(p,q)$,  for every $0\leq
t \leq d(p,q)$, we have that $d(\ga_{p,q}(t), \ga_{p,y}(t)) \leq
\delta$. (This actually holds for all $0\leq t\leq (y\cdot
\ga_{p,q}(\infty))_p$.) Now
\begin{align*}
&\abs{d^{\Phi}(\ga_{p,y}((x\cdot y)_p),\ga_{p,y}(d(p,q))))
-d^{\Phi}( \ga_{p, q}((x\cdot q)_p),q)} \leq D(3\delta).\end{align*}
This finishes the proof. \end{proof}

\begin{corollary}\label{cor:spikes} Assume $\Phi\in\mathcal{H}$ has $0$ pressure and
$G<\Isom(H)$ acts cocompactly. Fix $p \in H$. Then there exists a
constant $C_{\Phi}$ such that for all
$q \in H$, the $5$-tuple.
$$\left(\frac{d\mu_q^{\Phi}}{d\mu_p^{\Phi}},\ga_{p,q}(\infty), e^{-
d(p,q)}, d(p,q),C_{\Phi}\right)$$ is a ${\mathsf
G}^{Sym(\Phi)}$-spike.
\end{corollary}

\begin{proof} By assumption $\la_{\Phi}=0$, so by Proposition
\ref{prop:positive_average}, $\Phi$ has positive geodesic average at
all $p$. Multiply the ${\mathsf G}^{Sym(\Phi)}$-spike
$(e^{-2d^{Sym(\Phi)}( \ga_{p, q}((\xi\cdot q)_p),q)},
\ga_{p,q}(\infty), e^{- d(p,q)}, d(p,q),C_{Sym(\Phi)})$ by
$e^{d^{\Phi}(p,\ga_{p, q}((\xi\cdot q)_p))}$. By Lemma
\ref{lem:comparable_spikes}
$$ (e^{d^{\Phi}(p,\ga_{p,q}(\xi\cdot q)_p)-d^{\Phi}(q,\ga_{p,q}(\xi\cdot
q)_p)},\ga_{p,q}(\infty), e^{- d(p,q)}, d(p,q),C_{Sym(\Phi)})$$ is a
${\mathsf G}^{Sym(\Phi)}$-spike. By Lemma \ref{lem:approx_busemann}
and Lemma \ref{lem:comparable_spikes}
$$
\left(\frac{d\mu^{\Phi}_q}{d\mu_p^{\Phi}}(\xi)=e^{-\rho^{\Phi}_{\xi}(p,q)},\ga_{p,q}(\infty),
e^{- d(p,q)}, d(p,q),e^{6D(\delta)}C_{Sym(\Phi)}\right)$$ is a
${\mathsf G}^{Sym(\Phi)}$-spike. Set
$C_{\Phi}=e^{6D(\delta)}C_{Sym(\Phi)}$.
\end{proof}
\newpage
\section{Basis Theorem.}
In this section we prove an approximation theorem about ${\mathsf
G}$-spikes. On first glance, it may seem disconnected from our
stated goal in the introduction, but it plays the lead role in the
whole construction.

\begin{theorem} \label{thm:basis}
Let $(X,d,\nu)$ be a probability metric space.  Assume ${\mathsf
G}:X\times X \times \R_{\geq 0} \to \R_{\geq 0}$ is an almost
decreasing nicely decaying function with respect to $\nu$ with
constants $C_{\mathsf G}$, $\alpha_{\mathsf G}, \beta_{\mathsf G}$.
Assume
$\mathcal{F}=\{(f_{\alpha},r_{\alpha},a_{\alpha},s_{\alpha},C_{\alpha})\}_{\alpha\in
I}$ is a family of continuous unit ${\mathsf G}$-spikes.

Let $h:\R\to \R$ be some non-increasing function such that
$\lim_{t\to \infty}h(t)=0$. Assume that there exist constants $B \in
\nn$ and $\delta>0$ (we assume that $\delta$ and $R$ are the
constants for ${\mathsf G}$ from the definition of almost
decreasing) such that for every $D>0 $ and some $S\geq 0$ there
exists a countable subcollection $\{(f_i, r_i, a_i,s_i,
C_i)\}_{i=1}^{\infty}\subset \mathcal{F}$ such that
\begin{itemize}
\item{$r_i\leq h(S)$ and $s_i\geq S$ and $C_i\leq D$,}
\item{$\#\{i\,:\,x \in \Pi(a_i,r_i)\}\leq B$ for every $x \in
X$,}\item{$\abs{s_i-s_j}\leq \delta$ for all $i,j$,}
\item{$\nu(X-\bigcup_{i=1}^{\infty}\Pi(a_i,r_i)) \leq h(D),$}
\end{itemize}
\begin{enumerate}\item{ Then for every continuous uniformly positive function $F$ there exists a
countable set $\{\alpha_i\} \ssu I$ and constants
$\la_{\alpha_i}\geq 0$ such that
$$F=\sum_{i=1}^{\infty} \la_{\alpha_i} f_{\alpha_i},$$
with convergence in $L^1(X,\nu)$.} \item{Moreover, if there exists
one such $D>0$ and positive constants $m,\delta$ and $a$ such that
for every $S$ we can choose a finite subcollection $\{(f_i, r_i,
a_i,s_i, C_i)\}_{i=1}^{k}\subset\mathcal{F}$ with $C_i<D$ and the
extra conditions
\begin{itemize}
\item{ $m e^{-S} \leq r_i \leq e^{-S}$, and $S \leq s_i\leq
S+\delta$,}
\item{$X=\bigcup_{i=1}^{k}\Pi(a_i,r_i)$, and}
\item{$f_i$
are $a$-H\"{o}lder $\G$-spikes,} \end{itemize} then for every
positive H\"{o}lder function $F$ we can find a countable subset of
indices $\set{\alpha_i} \ssu I$ such that
$$F=\sum_{i=1}^{\infty} \la_{\alpha_i} f_{\alpha_i},$$
with convergence both pointwise and uniformly, and
$$\sum_{i=1}^{\infty} \la_{\alpha_i}\norm{f_{\alpha_i}}_{1,X}s_{\alpha_i} \leq \infty.$$}
\end{enumerate}
\end{theorem}
First, we will need to prove a Key Proposition.
\begin{proposition} \label{prop:subfunction} Let $(X,d,\nu)$ be a probability metric space.  Assume ${\mathsf G}:X\times
X \times \R_{\geq 0} \to \R_{\geq 0}$ is almost decreasing nicely
decaying with respect to $\nu$ function with constants $C_{\mathsf
G}$.
 For every
uniformly positive function $F$ we set
$$ t_{\infty}(F) = \sup\left\{\frac{F(y)}{F(x)}\,:\, x \in X, y \in Y\right\}$$ and
$$t_{\infty}(F) \geq t_{\ep} =\sup\left\{\frac{F(y)}{F(x)}\,:\, x \in X, y \in Y, \pi(x,y)\leq \ep\right\}\geq 1,$$
Choose $\ep>0$ and $S$ such that $t_{\infty}e^{-\beta_{\mathsf G}S}
\leq \ep^{\beta_{\mathsf G}} t_{\ep}$.

Assume that $\{(f_i, r_i, a_i,s_i, C_i)\}_{i=1}^{\infty}$ is a set
of $\G$-spikes such that
\begin{itemize}
\item{$r_i\leq \ep$  and $s_i\geq S$ and $C_i\leq D$,}
\item{$\#\{i\,:\,x \in \Pi(a_i,r_i)\}\leq B$ for every $x \in
X$,}\item{$\abs{s_i-s_j}\leq \delta$ for all $i,j$,}
\end{itemize}
and denote $Y=\bigcup_{i=1}^{\infty}\Pi(a_i,r_i)$.

 Then there exists constants $0\leq
\la_i=\frac{F(a_i)}{2DC_{\mathsf G} t_{\ep}^2B}$ such that
$$\sum_{i=1}^{\infty} \la_i f_i(x) \leq F(x) \text{ for all } x
\in X,$$ and
$$\sum_{i=1}^{\infty} \la_i f_i(x) \geq \frac{F(x)}{2D^2C_{\mathsf G} t_{\ep}^3B} \text{ for all } x
\in Y,$$
\end{proposition}

\medskip
\begin{proof}  Let $\hat{s} =\sup(s_i)\geq S$. Let us estimate
$\sum_{i=1}^{\infty} F(a_i) f(x)$.
\begin{align*}
\sum_{i=1}^{\infty} F(a_i) f_i(x) &\leq \sum_{i=1}^{\infty} F(a_i) D
e^{\alpha_{\mathsf G} s_i} \int_{\Pi(a_i,r_i)}{\mathsf {\mathsf
G}}^{\Phi}(x,y,s_i)d\nu(y)
\\ & \leq D\sum_{i=1}^{\infty} e^{\alpha_{\mathsf G} \hat{s}}
\int_{\Pi(a_i,r_i)}F(a_i)C_{\mathsf G}{\mathsf
G}^{\Phi}(x,y,\hat{s})d\nu(y) \\ & \leq DC_{\mathsf
G}\sum_{i=1}^{\infty} e^{\alpha_{\mathsf G} \hat{s}}
\int_{\Pi(a_i,r_i)}F(y){\mathsf G}^{\Phi}(x,y,\hat{s})d\nu(y)
\\& \leq DC_{\mathsf G} t_{\ep}Be^{\alpha_{\mathsf G}
\hat{s}} \int_{Y}F(y){\mathsf G}^{\Phi}(x,y,\hat{s})d\nu(y) \\
& \leq DC_{\mathsf G} t_{\ep}Be^{\alpha_{\mathsf G} \hat{s}}\left(
\int_{\Pi(x,\ep)}F(y){\mathsf G}^{\Phi}(x,y,\hat{s})d\nu(y)+\int_{Y- \Pi(x,\ep)}F(y){\mathsf G}^{\Phi}(x,y,\hat{s})d\nu(y)\right) \\
&\leq DC_{\mathsf G} t_{\ep}Be^{\alpha_{\mathsf G} \hat{s}}\left(
\int_{\Pi(x,\ep)}t_{\ep}F(x){\mathsf
G}^{\Phi}(x,y,\hat{s})d\nu(y)+\int_{Y-
\Pi(x,\ep)}F(x)t_{\infty}{\mathsf G}^{\Phi}(x,y,\hat{s})d\nu(y)\right)  \\
& \leq DC_{\mathsf G} t_{\ep}Be^{\alpha_{\mathsf G}
\hat{s}}F(x)\left(t_{\ep}e^{-\alpha_{\mathsf G}\hat{s}} + t_{\infty}
e^{-\alpha_{\mathsf G} \hat{s}} \frac{e^{-\beta_{\mathsf G}
\hat{s}}}{\ep^{\beta_{\mathsf G}}}\right) = 2DC_{\mathsf G}
t_{\ep}^2B F(x)
\end{align*}
Also observe that for all $y \in Y$ such that $\pi(x,a_j)\leq r_i
\leq \ep$ we have $$\sum_{i=1}^{\infty} F(a_i) f_i(x) \geq
F(a_j)\frac{f(a_j)}{D}\geq \frac{F(a_j)}{D} \geq
\frac{F(x)}{D\,t_{\delta}}.$$ Dividing by $2DC_{\mathsf G}
t_{\ep}^2B$ we obtain the proposition. \end{proof}

\medskip
Now we can prove Theorem \ref{thm:basis}.

\begin{proof} Let us describe the inductive procedure, how we
approximate $F$. Fix $\ell >1$.  Take a sequence $D_n\geq 1$ such
that $h(D_n) \to 0$ but $\sum_{n=1}^{\infty}\frac{1}{D_n^2} =
\infty$. Also fix any bounded sequences $\tau_n>1$ (from above) and
$\ga\leq\ga_n<1$ (we can easily assume that $\tau_n=2$ and $\ga_n=
1/2$).  Set $R_0=F$ and $\ep_{-1}=1$. Assume we have a continuous
function $R_n$ with respect to the distance $\pi$.

We set $g(x)=m\, x$. Construction steps: For every step we will have
$t_{\ep_n}\leq \ell$.
\begin{enumerate}
\item{ Find $\ep_{n}$ such that $t_{\ep_n}(R_n)\leq \ell$. In the
case $2$ i.e., when $R_n$ is $a$-H\"{o}lder set
$\ep_n^a=\inf\set{\frac{(\ell-1)\inf_{x \in X} R_{n}(x)}{\sup_{y \in
X}D_{g(\ep_{n-1})}^a R_n(y)}, g(\ep_{n-1})^a}$.} \item{Set $S_n$ to
be a number such that $h(S_n)\leq \ep_n$ and
$t_{\infty}(R_n)e^{-\beta_{\mathsf G}S_n}\leq \ep_n^{\beta_{\mathsf
G}}\ell. $ In case $h(s)\leq e^{-s}$ second inequality implies the
first one so we can actually set $e^{-S_n}=
\ep_n\left(\frac{\ell}{t_{\infty}(R_n)}\right)^{\frac{1}{\beta_{\mathsf
G}}}$}
\item{Find a set of spikes satisfying the conditions of Theorem
\ref{thm:basis} for $D_n$ and $S_n$,i.e, a set of spikes
$\{(f_{{\alpha_i}^{(n)}}, r_{{\alpha_i}^{(n)}},
a_{{\alpha_i}^{(n)}},s_{{\alpha_i}^{(n)}},
C_{{\alpha_i}^{(n)}})\}_{i=1}^{\infty}$ such that
\begin{itemize}
\item{$r_{\alpha_i^{(n)}}\leq h(S_n)$ and $s_{\alpha_i^{(n)}} \geq
S_n$ and $C_{\alpha_i^{(n)}}\leq D_n$,} \item{$\#\{i\,:\,x \in
\Pi(a_{\alpha_i^{(n)}},r_{\alpha_i^{(n)}})\}\leq B$ for every $x \in
X$,}\item{$\abs{s_{\alpha_i^{(n)}}-s_{\alpha_j^{(n)}}}\leq \delta$
for all $i,j$,( in the second case $S_n\leq s_i \leq S_n+\delta$)}
\item{$\nu(X-\bigcup_{i=1}^{\infty}\Pi(a_{\alpha_i^{(n)}},r_{\alpha_i^{(n)}}))
\leq h(D_n),$}
\end{itemize}}
\item{Find a finite a number $k_n$ such that
$\nu(X-\bigcup_{i=1}^{k_n}\Pi(a_{\alpha_i^{(n)}},r_{\alpha_i^{(n)}}))
\leq \tau_n h(D_n)$, and set
$Y_n=\bigcup_{i=1}^{k_n}\Pi(a_{\alpha_i^{(n)}},r_{\alpha_i^{(n)}})$.
In the second case this is automatic and $Y_n=X$.}
 \item{ By Proposition \ref{prop:subfunction} there
exists a function $h_n(x) = \sum_{i=1}^{k_n} \la_{\alpha_i^{(n)}}
f_{\alpha_i^{(n)}}(x)$ such that $\la_i> 0$ and \begin{align*}
h_n(x)&\leq F(x) &\forall x
\in X \\
h_n(x)&\geq \frac{F(x)}{2D_n^2 \ell^3 C_{\mathsf G}B} &\forall y \in
Y_n
\end{align*}}
\item{Set $R_{n+1} = R_n - \ga_n h_n$. From step $5$ and since
$\ga_n<1$, we have $R_{n+1}$ is uniformly positive and continuous.}
\end{enumerate}

 Now we estimate the $L^1$-norm of $R_{n+1}$ (Here $\norm{.}_{p,A}$ defines $L^p(A,\nu)$-norm for
 $1\leq p \leq \infty$):
 \begin{align*}
 \norm{R_{n+1}}_{1,X} &= \norm{R_n-\ga_n h_n}_{1,X} \leq \norm{R_n-\ga_n
 h_n}_{1,Y_n}+\norm{R_n-\ga_n h_n}_{1,X-Y_n}  \\ & \leq
 \left(1- \frac{\ga_n}{2D_n^2 \ell^3 C_{\mathsf G}B} \right)\norm{R_n}_{1,Y_n}+\nu(X-Y_n)\norm{R_n}_{\infty,X-Y_n}
\\ & \leq \left(1- \frac{\ga_n}{2D_n^2 \ell^3 C_{\mathsf G}B} \right)\norm{R_n}_{1,X}+\tau_n
h(D_n)\norm{F}_{\infty,X}
 \end{align*}
Since $\tau_n h(D_n) \to 0$ and $\prod_{n=1}^{\infty}\left(1-
\frac{\ga_n}{2D_n^2 \ell^3 C_{\mathsf G}B} \right) =0$  we obtain
(by Lemma 5.12 of \cite{ConnellMuchnik04}, for instance) that
$\norm{R_n}_{L^1(X,\nu)}\to 0$. From construction we have
$$F(x)=R_0(x) = \sum_{n=0}^{\infty} \ga_n h_n(x) =
\sum_{n=0}^{\infty} \ga_n \sum_{i=1}^{k_n} \la_{\alpha_i^{(n)}}
f_{\alpha_i^{(n)}}(x) \text{ in } L^1(X,\nu).$$

In the second $2$ we have
\begin{align*}
R_{n+1}(x)& = R_n(x) -\ga_n h_n(x) \leq \left(1- \frac{\ga_n}{2D_n^2
\ell^3 C_{\mathsf G}B} \right)R_n(x) \\ & \leq
\prod_{i=0}^{n}\left(1- \frac{\ga_i}{2D_i^2 \ell^3 C_{\mathsf G}B}
\right) F(x). \end{align*} This proves that
$$F(x)= \sum_{n=1}^{\infty} \ga_n
h_n(x) = \sum_{n=1}^{\infty} \ga_n \sum_{i=1}^{k_n}
\la_{\alpha_i^{(n)}} f_{\alpha_i^{(n)}}(x) \text{ with uniform
convergence.}$$

Now let us prove the second estimate in case $2$.  First few of
trivial observations: \begin{itemize} \item{
$$\sum_{i=1}^{k_n} \la_{\alpha_i^{(n)}}\norm{f_{\alpha_i^{(n)}}}_{1,X}=\norm{h_n}_{1,X} \leq \norm{R_n}_{1,X} \leq \prod_{i=0}^{n-1}\left(1- \frac{\ga_i}{2D_i^2 \ell^3 C_{\mathsf G}B} \right)
\norm{F(x)}_{1,X}.$$} \item{ \begin{align*}t_{\infty}(R_{n+1})
&=\sup_{x,y\in X} \frac{R_{n+1}(x)}{R_{n+1}(y)} = \sup_{x,y\in X}
\frac{R_{n}(x)-\ga_nh_n(x)}{R_{n}(y)-\ga_nh_n(y)}\\ & \leq
\sup_{x,y\in X} \frac{\left(1-\frac{\ga_n}{2D_i^2 \ell^3 C_{\mathsf
G}B} \right)R_{n}(x)}{(1-\ga_n)R_{n}(y)} =
\frac{\left(1-\frac{\ga_n}{2D_i^2 \ell^3 C_{\mathsf G}B}
\right)}{1-\ga_n}t_{\infty}(R_n). \end{align*}} \item{ For all $x,y
\in X$ such that $\pi(x,y)\leq \ep_n$ with $\ep_n$ as defined we
have
\begin{align*}
t_{\ep_n}(R_n)=&\sup_{\substack{ y\not=x \\ \pi(x,y)\leq
\ep_n}}\frac{R_n(x)}{R_n(y)}=\sup_{\substack{ y\not=x \\
\pi(x,y)\leq \ep_n}}\frac{\abs{R_n(x)-R_n(y)}}{R_n(y)}+1 \\ & \leq
\sup_{\substack{ y\not=x \\
\pi(x,y)\leq \ep_n}}\frac{\abs{R_n(x)-R_n(y)}}{\pi(y,x)^a}
\cdot\frac{\ep^a}{\inf_{y \in X}R_n(y)}+1 \\ & \leq D_{\ep_n}^a(R_n)
\frac{\ell-1}{\sup_{y \in X} D_{\ep_{n-1}}^a R_{n}(y)}+1 \leq \ell
\intertext{since $\ep_n \leq \ep_{n-1}$.}
\end{align*}}
\end{itemize}

In view of this we will estimate $D_{g(\ep_n)}^a R_{n+1}(x)$.
\begin{align*}
D_{g(\ep_n)}^ah_n(x) \leq \sum_{i=1}^{k_n}
\la_{\alpha_i^{(n)}}D_{g(\ep_n)}^a f_{\alpha_i^{(n)}}(x) \leq
\sum_{i=1}^{k_n} \la_{\alpha_i^{(n)}}D_n\frac{
f_{\alpha_i^{(n)}}(x)}{g(\ep_n)^a} = D_n\frac{ h_n(x)}{g(\ep_n)^a}.
\end{align*}
Since $g(\ep_n)\leq \ep_n \leq g(\ep_{n-1})$ we have
\begin{align*} D_{g(\ep_n)}^a R_{n+1}(x) &= D_{g(\ep_n)}^a
R_{n}(x)+\ga_n D_{g(\ep_n)}^ah_n(x) \leq D_{g(\ep_{n-1})}^a
R_{n}(x)+\ga_n D_n\frac{h_n(x)}{g(\ep_n)^a} \\ & \leq
\frac{(\ell-1)inf_{x \in X} R_{n}(x)}{\ep_n^a} + \ga_n D_n\frac{
h_n(x)}{g(\ep_n)^a}  \\ & \leq
\max(\ell-1,D_n)\left(\frac{R_{n}(x)}{\ep_n^a} + \ga_n \frac{
h_n(x)}{g(\ep_n)^a}\right) =
(\ell+D_n)\frac{R_{n+1}(x)}{g(\ep_{n})^a}
\end{align*}

Now \begin{align*} \ep_{n+1}^a &=\inf\left(\frac{(\ell-1)inf_{x \in
X} R_{n+1}(x)}{sup_{y \in X}D_{g(\ep_{n})}^a R_{n+1}(y)},
g(\ep_{n})^a\right) \\ & \geq\inf\left(\frac{(\ell-1)inf_{x \in X}
R_{n+1}(x)}{(\ell+D_n) sup_{x\in X}
R_{n+1}(x)}g(\ep_n)^a, g(\ep_{n})^a\right) \\
& =\frac{(\ell-1)inf_{x \in X} R_{n+1}(x)}{(\ell+D_n) sup_{x\in X}
R_{n+1}(x)}g(\ep_n)^a = \frac{\ell-1}{(\ell+D_n)}
\frac{g(\ep_n)^a}{t_{\infty}(R_{n+1})} \\ & \geq
\frac{(\ell-1)m}{(\ell-1+D_n)} \frac{\ep_n^a}{t_{\infty}(R_{n+1})}.
\end{align*}
Simple induction proves that \begin{align*}\ep_{n+1}^a \geq
\left(\frac{(\ell-1)m}{(\ell+D_n)}\right)^n
\frac{\ep_0^a}{\prod_{i=1}^{n+1}t_{\infty}(R_i)}
\end{align*}
Recall that $D_n \leq D$ and we set $\ga_n=\ga$ to be a number such
that $\frac{\left(1-\frac{\ga_n}{2D_i^2 \ell^3 C_{\mathsf G}B}
\right)}{1-\ga_n} = K\geq 1$.  Thus we have $t_{\infty}(R_n)\leq K^n
t_{\infty}(R_0)$. And so
$$\ep_{n+1}^a\geq \left(\frac{(\ell-1)m}{(\ell+D_n)}\right)^n
\frac{\ep_0^a}{\prod_{i=1}^{n+1}K^it_{\infty}(R_0)} \geq
\left(\frac{(\ell-1)m}{(\ell+D)}\right)^n\frac{K^{-\frac{(n+1)(n+2)}{2}}}{t_{\infty}(R_0)^{n+1}}.$$

Now
\begin{align*}
e^{-aS_n} =
e^{a}\left(\frac{\ell}{t_{\infty}(R_n)}\right)^{\frac{a}{\beta_{\mathsf
G}}} \geq
\left(\frac{(\ell-1)m}{(\ell+D)}\right)^n\frac{K^{-\frac{(n+1)(n+2)}{2}}}{t_{\infty}(R_0)^{n+1}}
\left(\frac{\ell}{K^n
t_{\infty}(R_0)}\right)^{\frac{a}{\beta_{\mathsf G}}}
\end{align*}
 Now we have \begin{align*} s_{\alpha_i^{(n)}} \leq S_n+\delta
\leq p(n),
\end{align*}
where $p(x)$ is a polynomial of degree $2$. Now \begin{align*}
\sum_{n=0}^{\infty} \sum_{i=1}^{k_n} \la_{\alpha_i^{(n)}}
\norm{f_{\alpha_i^{(n)}}}_{1,X}s_{\alpha_i^{(n)}} \leq
\sum_{n=0}^{\infty} (S_n+\delta) \left(1-\frac{\ga}{2D^2 \ell^3
C_{\mathsf G}B} \right)^n \norm{R_0}_{1,X} <\infty,
\end{align*}
as $S_n$ grow polynomially. \end{proof}

\begin{proof}[Proof of Theorem \ref{thm:main2}]
We wish to apply the above theorem by setting $F=f$, $X=\pa H$,
$\nu=\nu_p^{\Phi}$, and the index set $I$ to be the group $G$. Under
the hypotheses, the radon Nikodym derivatives, $\hat{f}_g=g_*F
\frac{dg_*\nu}{d\nu}$ are $\mathsf{G}^{\op{Sym}(\Phi)}$-spikes by
Corollary \ref{cor:spikes} and Lemma \ref{lem:comparable_spikes}. We
make them unit spikes by dividing by
$\hat{f}_g(\ga_{p,gp}(\infty))=F(\ga_{g^{-1}p,p}(\infty))e^{d_{\Phi}(p,gp)}$
to obtain unit spikes $f_g$. Under the assumptions of the theorem,
$H$ has bounded flip. Therefore, we may apply Proposition
\ref{prop:RN_holder} to show that $\frac{dg_*\nu}{d\nu}$ is
H\"{o}lder with $D^\eps_{r_g} \frac{dg_*\nu}{d\nu}\leq
\frac{C}{r_g^\eps}$. On the other hand, the action of $g$ is
conformal and Lipschitz with respect to $\pi_p$ with $D_r
g(\zeta)=\sup_{\xi\in
\Pi(\zeta,r)\setminus\zeta}e^{-B_\xi(p,gp)}\leq e^{d(p,gp)}$. Hence
by Lemma \ref{lem:Holder_const}, applying the chain rule and product
rule we have that $D_{r_g}^\eps f_g\leq \frac{C}{r_g^\eps}$ with $C$
independent of $g$. Hence the $f_g$ form a family of
$\eps$-H\"{o}lder $\G$-spikes. Moreover, the spike constants are
related by $r_g=e^{-s_g}=e^{-d(p,gp)}$ and hence for any $S>0$ and
$\delta>\op{diam}(H/G)$ the family of shadows $\Pi(a_g,r_g)$ for
$S\leq s_g\leq S+\delta$ cover all of $\pa H$.

Therefore we obtain both conclusions of Theorem \ref{thm:basis}.
Setting $\mu_p(g)=\la_g \norm{f_g}_1=\la_g e^{-d_{\Phi}(p,gp)}$, the
first statement of the theorem states that $\mu\conv f\nu=f\nu$. The
second part states that first moment is finite, $\sum_{g\in G}
\mu_p(g)\log d(p,gp)<\infty$. We may write the entropy of $\mu_p$ as
$-\sum_{g\in G} \mu_p(g)\log \mu_p(g)=\sum_{g\in G}
\mu_p(g)d_{\Phi}(p,gp)-\sum_{g\in G} \la_ge^{-d_{\Phi}(p,gp)} \log
\la_g$. Since $d_{\Phi}\leq \norm{\Phi}_\infty d$, this is finite if
and only if $\sum_{g\in G} \la_ge^{-d_{\Phi}(p,gp)} \log
\frac{1}{\la_g}<\infty$. On the other hand, $\sum_{g\in G}
\la_ge^{-d_{\Phi}(p,gp)} d(p,gp)<\infty$. So the former sum
converges unless possibly for some subsequence $g_i\in G$
$\lim_{i\to\infty} \frac{-\log \la_{g_i}}{d(p,g_ip)}\to \infty$.
However, for any such subsequence the sum $-\sum_i \la_{g_i} \log
\la_{g_i}$ is bounded since the number of elements with $d(p,gp)$
fixed grows at most exponentially while $\la_{g_i}$ decays
superexponentially in $d(p,gp)$. Hence the entropy of $\mu_p$ is
finite. By the criterion of Kaimanovich in \cite{Kaimanovich00}, the
measure $f\nu^\Phi_p$ is a Poisson boundary for $\mu_p$.
\end{proof}

\begin{remark}
By following the procedure found in \cite{ConnellMuchnik04}, we can
vary the $\la_g$ at each stage in the construction in Theorem
\ref{thm:basis}. This allows us to produce an infinite dimensional
space of measures $\mu$ (random walks) on $G$ each with the same
Poisson Boundary $(\pa H,f\nu_p^\Phi)$.
\end{remark}

\newpage
\appendix
\section{Unit Tangent Spaces}\label{sec:appendixA}



We recall the notations from Section \ref{sec:Gibbs_States}.

\begin{lemma}\label{lem:dist_form} Assume $H$ is $CAT(-1)$-space. Let $p \in H$ and $\ga_1(t)$ and $\ga_2(t)$ be two
geodesics such that $\ga_1(0)=\ga_2(0)=p$.  Then for all $0 \leq s
\leq S$ and $0\leq t \leq T$ we have,

\begin{align*}
\frac{\cosh\(d(\ga_1(s),
\ga_2(t))\)-\cosh(t-s)}{\cosh(t)\cosh(s)-\cosh(t-s)} &\leq
\frac{\cosh\(d(\ga_1(S),
\ga_2(T))\)-\cosh(T-S)}{\cosh(T)\cosh(S)-\cosh(T-S)},
\end{align*}
and if $s=t$ and $S=T$ we have
\begin{align*}
d(\ga_1(t), \ga_2(t)) &\leq
2\op{arcsinh}\(\sinh\left(\frac{d(\ga_1(T),
\ga_2(T))}{2}\right)\frac{\sinh(t)}{\sinh(T)}\).
\end{align*}
In particular,  $d(\ga_1(t), \ga_2(t))\leq 2
\sinh\left(\frac{d(\ga_1(T), \ga_2(T))}{2}\right)e^{t-T}.$ This is
accurate when $T \gg t+\frac12 d(\ga_1(T), \ga_2(T))$. Similarly,
when $t$ and $d(\ga_1(T), \ga_2(T))$ are large, we have the estimate
$\frac{d(\ga_1(T),\ga_2(T))-d(\ga_1(t), \ga_2(t))}{T-t}\geq
2-O\(\max\set{e^{-t},e^{-d(\ga_1(T),\ga_2(T))}}\).$
\end{lemma}

\medskip
\begin{proof}
1)We use the law of cosines for the hyperbolic plane of curvature
$-1$: for a triangle with sides $a,b,c$ we have

\begin{align*}
\cosh(c)&=\cosh(a)\cosh(b) -
\sinh(a)\sinh(b)\cos(\theta)\\
&=\cosh(a-b)+\sinh(a)\sinh(b)\(1-\cos(\theta)\),
\end{align*}
where $\theta$ is an angle between sides $a$ and $b$.

Now suppose $\sigma_1$ and $\sigma_2$ are unit speed geodesics in
$\HH^2$ with $\sigma_1(0)=\sigma_2(0)$ and
$d_{\HH^2}\(\sigma_1(S),\sigma_2(T)\)=d(\ga_1(S),\ga_2(T)).$ The
triangles $[\sigma_1(0),\sigma_1(s),\sigma_2(t)]$ and
$[\sigma_1(0),\sigma_1(S),\sigma_2(T)]$ in $\HH^2$ share the same
angle at $\sigma_1(0)$. Therefore, if $a=s,b=t$ and
$c(s,t)=d_{\HH^2}\(\sigma_1(s), \sigma_2(t)\)$ in the law of
cosines, then we have
$$\frac{\cosh(c(s,t))-\cosh(s)\cosh(t)}{\cosh(c(S,T))-\cosh(S)\cosh(T)}=
\frac{\sinh(s)\sinh(t)}{\sinh(S)\sinh(T)}=
\frac{\cosh(t)\cosh(s)-\cosh(t-s)}{\cosh(T)\cosh(S)-\cosh(T-S)}.$$

By the comparison property for the CAT(-1) space $H$ for a triangle
with side lengths $S,T,$ and $d\(\ga_1(S), \ga_2(T)\)$, we have that
$d\(\ga_1(s), \ga_2(t)\)$ is not more than the same distance
measured in the hyperbolic triangle with the same side lengths.
Since $\cosh$ is increasing we obtain the first estimate. For the
second we take $s=t$ and $S=T$ to obtain,
$$\cosh\(d\(\ga_1(t), \ga_2(t)\)\)-1 \leq (\cosh(d(\ga_1(T),
\ga_2(T)))-1) \left(\frac{\sinh(t)}{\sinh(T)}\right)^2.$$

Since $\cosh(x)-1 = 2 \sinh(x/2)^2$ we obtain the second estimate.
Using $\sinh(x)\geq x$ for $x\geq 0$, this implies
$$\frac{d\(\ga_1(t), \ga_2(t)\)^2}{2} \leq 2 \sinh\left(\frac{d(\ga_1(t), \ga_2(t))}{2}\right)^2
\leq 2\sinh\left(\frac{d(\ga_1(T),
\ga_2(T))}{2}\right)^2\left(\frac{\sinh(t)}{\sinh(T)}\right)^2. $$
Taking the square root and using the remark gives the second
estimate. The last statement follows by expanding the sharp estimate
by a power series in $e^{-t}$ and $e^{- d(\ga_1(T),\ga_2(T))}$ and
using the fact that $2T\geq d(\ga_1(T),\ga_2(T))$ and $T\geq t.$
\end{proof}


We will sometimes distinguish a geodesic in $SH$ which starts at
$x\in H$ by the notation $\ga_x$ and a geodesic starting at $x$ and
passing through $y\in \bar{H}$ by $\ga_{x,y}$. One of the principal
differences between CAT$(-1)$ spaces and negatively curved manifolds
is that while in both cases the geodesic segment between $x$ and $y$
is unique, in the former case there may be an infinite number of
geodesics passing through any two points.

\begin{corollary}\label{cor:horodist} Assume that $H$ is $CAT(-1)$-space.  Let $\ga_1$
and $\ga_2$ be two geodesics such that $\ga_1(\infty)=
\ga_2(\infty)=\zeta \in\pa H$. Setting $d=d(\ga_1(0),\ga_2(0))$ and
$\rho=\rho^0_{\zeta}(\ga_1(0),\ga_2(0))$, we have
$$d(\ga_1(s),\ga_2(t))\leq
\cosh^{-1}\(\frac{\cosh(d)-\cosh(\rho)}{e^{t+s}}+\cosh(\rho+s-t)\).$$

Assume $\rho^0_{\zeta}(\ga_1(0),\ga_2(0))=0$ and $s=t$, Then
$$d(\ga_1(t),\ga_2(t)) \leq 2 \sinh^{-1}\(\sinh(d/2) e^{-t}\)\leq
\begin{cases}
d-\frac{2}{d}\(e^{-d}+d-1\)t & 0\leq t\leq \frac{d}{2}\\
2\sinh\left(\frac{d}{2}\right)e^{-t}& t>\frac{d}{2}\end{cases}.$$
\end{corollary}

\medskip
\begin{proof} Let $\ga_1(0)=p, \ga_2(0)=q$. Choose a point $z_T$ such
that $d(p,z_T)=d(q,z_T)=T$.  Then
\begin{align*}
d(\ga_{p,z_T}(t),\ga_{q,z_T}(s))=d(\ga_{z_T,p}(T+c-t),\ga_{z_T,q}(T-s))
\end{align*}
Note that $z_T \to \zeta$ as $T \to \infty$ and so $\ga_{p,z_T}(t)$
and  $\ga_{q,z_T}(t)$ tend to $\ga_{p,\zeta}(t)$ and
$\ga_{q,\zeta}(t)$ respectively. Using the first estimate of the
previous lemma, and noting that
$\frac{\sinh(T+\rho-t)\sinh(T-s)}{\sinh(T+\rho)\sinh(T)}$ tends to
$e^{-t-s}$ as $S,T\to \infty$, we obtain the first inequality.

For the second inequality, take $\rho=0$ and $s=t$. In this case,
the distance reduces (after applying identities) to $2
\sinh^{-1}\(\sinh(d/2) e^{-t}\)$. Note that this function has
positive second derivative in $t$. Therefore by convexity, the line
between the initial value at $t=0$ and the upper estimate $2
\sinh(\frac{d}{2})e^{-t}$ at $t=\frac{d}{2}$ is an upper bound.
\end{proof}

\begin{lemma} Assume that $H$ is $CAT(-1)$-space.  Let $\ga_1$ and
$\ga_2$ be two geodesics with $\ga_1(0)= \ga_2(0)=p$.  If
$c_+=(\ga_1(\infty)\cdot \ga_2(\infty))_p$, then for any $s\geq 0$
we have
\begin{align*}
d(\ga_1(s),\ga_2(t))\leq
\cosh^{-1}\(e^{-c_+}\sinh(s)\sinh(t)+\cosh(t-s)\),
\end{align*}
and if $s=t$, then
\begin{align*}
d(\ga_1(t),\ga_2(t))\leq
2\sinh^{-1}\(e^{-c_+}\sinh(t)\)\leq\begin{cases} 2 e^{-c_+}\sinh(t)
& t\leq
c_+ \\
2(t-c_+)+2 e^{-t}\sinh(c_+) & t>c_+
\end{cases}.
\end{align*}
\end{lemma}

\begin{proof}
First observe that in the formula from Lemma \ref{lem:dist_form}, we
may replace $\frac12 d(\ga_1(T),\ga_2(T))$ by $T-c_+$ as we take
$S=T\to \infty$. Taking this limit we obtain
$$\frac{\cosh\(d(\ga_1(S),
\ga_2(T))\)-\cosh(T-S)}{\cosh(T)\cosh(S)-\cosh(T-S)}\to e^{-c_+}.$$
Then rearranging terms on the left hand side we obtain the first
inequality. Setting $s=t$ and using standard identities we obtain
$$d(\ga_1(s),\ga_2(t))\leq 2 \sinh^{-1}(e^{-c_+}
\sinh(t)).$$

Since $\sinh^{-1}(t)=t-\frac{t^3}{6}+O(t^4)$ and $\sinh^{-1}(1)<1$,
we have $\sinh^{-1}(t)\leq t$ for $0\leq t\leq 1$. So
$d(\ga_1(t),\ga_2(t))\leq 2 e^{-c_+}\sinh(t)$ for $0\leq t\leq
\sinh^{-1}(e^{c_+})$. On the other hand, setting
$f(t,c_+)=2(t-c_+)+2 e^{-t}\sinh(c_+)-2 \sinh^{-1}(e^{-c_+}
\sinh(t))$, we have
$$\frac{df(t,c_+)}{dc_+}=2\(-1+e^{-t}\cosh(c_+)+\(1+\frac{e^{2c_+}}{\sinh(t)^2}\)^{-1/2}\)$$
which is strictly positive. Since $f(t,0)=0$, we have
$$d(\ga_1(t),\ga_2(t))\leq 2(t-c_+)+2 e^{-t}\sinh(c_+)$$ for all
$t,c_+\geq 0$. However, the former estimate is more accurate when
$t<c_+$.
\end{proof}

\begin{lemma} Assume that $H$ is $CAT(-1)$-space.  Let $\ga_1$ and
$\ga_2$ be two geodesics with $\ga_1(0)= \ga_2(0)=p$.  If
$c_+=(\ga_1(\infty)\cdot \ga_2(\infty))_p$, then for any $s\geq 0$
we have
\begin{align*}
\int_{0}^{\infty} d(\ga_1(t),\ga_2(t))e^{-\abs{t-s}}dt \leq
4\max\set{s-c_+,0}
+e^{-\abs{c_+-s}}\(\abs{c_+-s}+3\)-(s+1)e^{-c_+-s}.
\end{align*}
\end{lemma}

\begin{proof}

Now assume $s\geq c_+$, then we have
\begin{align*}
&\int_{0}^{\infty} d(\ga_1(t),\ga_2(t))e^{-\abs{t-s}}dt  \\& =
\int_{0}^{c_+} d(\ga_1(t),\ga_2(t))e^{t-s}dt+\int_{c_+}^{s}
d(\ga_1(t),\ga_2(t))e^{-t+s}dt+\int_{s}^{\infty}
d(\ga_1(t),\ga_2(t))e^{-t+s}dt\\ \intertext{Using the estimate of
the previous lemma we can integrate each of these to obtain}
&\int_{0}^{\infty} d(\ga_1(t),\ga_2(t))e^{-\abs{t-s}}dt\leq
e^{s-c_+}(s-c_++3)-(s+1)e^{-c_+-s}+4(s-c_+).
\end{align*}
For $s\leq c_+$ we obtain,
\begin{align*}
&\int_{0}^{\infty} d(\ga_1(t),\ga_2(t))e^{-\abs{t-s}}dt  \\& =
\int_{0}^{s} d(\ga_1(t),\ga_2(t))e^{t-s}dt+\int_{s}^{c_+}
d(\ga_1(t),\ga_2(t))e^{-t+s}dt+\int_{c_+}^{\infty}
d(\ga_1(t),\ga_2(t))e^{-t+s}dt\\
&\leq e^{c_+-s}(c_+-s+3)-(s+1)e^{-c_+-s}.
\end{align*}

Combining both of these cases together completes the proof.
\end{proof}

%
%
\medskip

%

\begin{lemma} \label{lem:Sasaki_dist_estimate} Let $H$ be $CAT(-1)$-space. Let $\ga_1$ and $\ga_2$ are two geodesics, such that
$\ga_1(0)= \ga_2(0)=p$.  Let $c_+=(\ga_1(\infty)\cdot
\ga_2(\infty))_p$ and $c_-=(\ga_1(-\infty)\cdot \ga_2(-\infty))_p$.
Then for all $s>0$ we have
$$\op{dist}({\mathsf g}^s(\ga_1), {\mathsf g}^s(\ga_2) \leq
2\max\set{s-c_+,0} +\frac{\(\abs{c_+-s}+3\)}{2
e^{\abs{c_+-s}}}-\frac{(s+1)}{2e^{s+c_+}}+\frac{(c_-+2)}{2e^{s+c_-}}.
$$
\end{lemma}

\begin{proof}
\begin{align*}
\op{dist}({\mathsf g}^s\ga_1, {\mathsf
g}^s\ga_2)=&\frac{1}2\int_{-\infty}^{\infty} d(\ga_1(t+s),
\ga_2(t+s))e^{-|t|}dt =\frac{1}2
\int_{-\infty}^{\infty} d(\ga_1(t), \ga_2(t))e^{-\abs{t-s}}dt \\
=&\frac{1}{2}\int_{-\infty}^{0} d(\ga_1(t), \ga_2(t))e^{t-s}dt
+\frac{1}{2}\int_{0}^{\infty} d(\ga_1(t), \ga_2(t))e^{-\abs{t-s}}dt
\end{align*}

Now for $t\leq 0$ we use simple substitution $t=-t$ to get
\begin{align*}
&\int_{-\infty}^{0} d(\ga_1(t), \ga_2(t)) e^{t-s} dt \leq
e^{-s}\int_{0}^{\infty}d(\ga_1(-t), \ga_2(-t)) e^{-t} dt
\intertext{Now using the previous lemma,} &\leq
e^{-s}e^{-c_-}(c_-+2).
\end{align*}

Again by the previous lemma,
$$\int_{0}^{\infty} d(\ga_1(t), \ga_2(t))e^{-\abs{t-s}}dt \leq
4\max\set{s-c_+,0}
+e^{-\abs{c_+-s}}\(\abs{c_+-s}+3\)-(s+1)e^{-c_+-s}.$$
\end{proof}
%
%
%
%

\begin{lemma}\label{lem:integral_estimate} Let $H$ be $CAT(-1)$-space. Let $\ga_1$ and $\ga_2$ are two geodesics, such that
$\ga_1(0)= \ga_2(0)=p$. Let $c_+=(\ga_1(\infty)\cdot
\ga_2(\infty))_p$,  $c_-=(\ga_1(-\infty)\cdot \ga_2(-\infty))_p$ and
$0<\beta\leq 1$. Then for all $T>0$ we have
\begin{align*}
 \int_{0}^{T}\op{dist}({\mathsf g}^s\ga_1, {\mathsf g}^s\ga_2)^{\beta}ds \leq &
 \(\frac{1}{\beta}+\frac{c_-}{2}\)\(e^{-\beta\, c_-}-e^{-\beta
 T-\beta\,c_-}\)+\frac{2}{\beta}\(1+\op{sgn}(T-c_+)\(1-e^{-\beta\abs{T-c_+}}\)\)\\
 &-\(\frac{2}{\beta}+\frac{c_+}{2}\)e^{-\beta
 \,c_+}-\frac{T-c_+}{2}
 e^{-\beta \abs{T-c_+}} +\frac{2^\beta\,\max\set{\left(T -c_{+}
 \right),0}^{1 + \beta}}{1 + \beta}\\
 \leq & \(\frac{1}{\beta}+\frac{c_-}{2}\)e^{-\beta\,
 c_-}+\frac{2}{\beta}\min\set{2,e^{\beta (T-c_+)}}-\(\frac{2}{\beta}+\frac{c_+}{2}\)e^{-\beta
 \,c_+} \\
 &-\frac{T-c_+}{2}
 e^{-\beta \abs{T-c_+}} +\frac{2^\beta\,\max\set{\left(T -c_{+}
 \right),0}^{1 + \beta}}{1 + \beta}\\
         \leq & \frac{5}{\beta}+
 \max\set{\left(T -c_{+} \right),0}^{1 + \beta}.
\end{align*}
\end{lemma}

\medskip
\begin{proof}
Since $0<\beta\leq 1$, we use $(x+y)^{\beta} \leq
x^{\beta}+y^{\beta}$ for $x,y >0$ and $(a+x)^\beta\leq a+\beta x$
for $a\geq 1$ to obtain

$$\op{dist}({\mathsf g}^s\ga_1, {\mathsf g}^s\ga_2)^\beta\leq  \frac{2 + \beta
c_{-}}{2\,e^{\beta\,\( c_{-} + s \) }} +\frac{\( 3 +\beta \abs{c_{+}
- s }\)}{2\,e^{\beta\,\abs{ s-c_{+}}
}}+2^\beta\max\set{|s-c_+|,0}^\beta.$$ Integrating we obtain the
first estimate.

For the second inequality, we use $1+\op{sgn}(x)(1-e^{-\abs{x}})\leq
\min\set{2,e^{x}}$. For the third inequality, we overestimate the
last term using that $2^\beta\leq 1+\beta$ for $\beta\in [0,1]$.
Note that $\(\frac{1}{\beta}+\frac{c_-}{2}\)\(1-e^{-\beta
T}\)e^{-\beta c_-}\leq \frac{1}{\beta}$. So the worst case occurs in
the first four terms by taking $T\to \infty$, in which case we
recover the given value.
\end{proof}

%
%
%

\medskip
%
%

In Particular, we have the following special case.
\begin{corollary}\label{cor:ac_estimate}
Let $\ga_1, \ga_2$ are $2$ geodesics in $H$, such that
$\ga_1(0)=\ga_2(0)=p$ and $c_{\pm}=(\ga_1(\pm\infty)\cdot
\ga_2(\pm\infty))_p$. Then for any $0\leq \alpha\leq 1$,

\begin{align*}
\int_0^{\alpha\, c_+} \op{dist}({\mathsf g}^s \ga_1, {\mathsf
g}^s\ga_2)^{\beta}ds \leq &
 \min\set{\frac{1}{\beta},\alpha c_+}\(\(1 + \frac{\beta c_-}{2}\)\,e^{-\beta\,c_-}\,
       +2e^{-\beta \,(1-\alpha)\,c_+}\)\\
       &-\frac{c_+}{2}e^{-\beta
       c_+} +\frac{(1-\alpha)c_+}{2}e^{-\beta(1-\alpha)c_+}.\\
\end{align*}
If in addition we assume $c_-\geq (1-\alpha)c_+$, then we have
\begin{align*}
\int_0^{\alpha\, c_+} \op{dist}({\mathsf g}^s \ga_1, {\mathsf
g}^s\ga_2)^{\beta}ds \leq &
e^{-\beta(1-\alpha)c_+}\(\frac{3}{\beta}+(1-\alpha)c_+\).
\end{align*}
\end{corollary}
\begin{proof}
We use the first estimate of the previous Lemma with $T=\alpha
\,c_+$. We then note that $1-e^{-\beta x}\leq \min\set{1,\beta x}$
to obtain the first estimate. Using $c_-=(1-\alpha)c_+$ and dropping
the $-\frac{c_+}{2}e^{-\beta c_+}$ term we obtian the second
estimate.
\end{proof}


We finish this section with a simple lemma that we will use later.
\begin{lemma} \label{lem:increasing_dist} Let $H$ be $CAT(-1)$ space.  Assume that
$\ga_1(0)=\ga_2(0)=p$. Then for all $a\geq 0$ and $t\geq 0$ we have
$$d(\ga_1(a),\ga_2(a)) \leq d(\ga_1(a),\ga_2(a+t)).$$
\end{lemma}

\medskip
\begin{proof} Using the comparession triangle we have
\begin{align*}
\cosh(d(\ga_1(a),\ga_2(a+t))) \geq
\cosh(a)\cosh(a+t)-\sinh(a)\sinh(a+t)cos(\angle_p(\ga_1,\ga_2))
\end{align*}
Differentiating with respect $t$ we have
\begin{align*}
\frac{d(\cosh(d(\ga_1(a),\ga_2(a+t))))}{dt} \geq
\cosh(a)\sinh(a+t)-\sinh(a)\cosh(a+t)cos(\angle_p(\ga_1,\ga_2))\geq
0,
\end{align*}
for $t\geq 0$ as $\frac{\cosh(x)}{\sinh(x)}$ is a decreasing
function. \end{proof}

\newpage
\section{$c_{\ga}^{\Phi}(x,y)$ is H\"{o}lder.}

We recall the notations from Section \ref{sec:Gibbs_States}. We let
$\Phi$ be a tempered H\"{o}lder function on $SH$ with H\"{o}lder
exponent $0<\beta\leq 1$ and global H\"{o}lder constant $K$.

\begin{lemma} \label{lem:dist_difference} Let $p \in H$. Assume
$\ga_1, \ga_2$ are two geodesics, such that $\ga_1(0)=\ga_2(0)=p$
and $c_+=\(\ga_1(\infty)\cdot \ga_2(\infty)\)_p$. Then for all
$0\leq \alpha\leq 1$ we have
$$\abs{d^{\Phi}(p,\ga_1(\alpha c_+))-d^{\Phi}(p,\ga_2(\alpha c_+))}
\leq K\,\(\(\frac{1}{\beta} + \frac{c_-}{2}\)\,e^{-\beta\,c_-}\,
       +\(\frac{2}{\beta}+\frac{(1-\alpha)c_+}{2}\)e^{-\beta \,(1-\alpha)\,c_+}\).$$
\end{lemma}

\medskip
\begin{proof}
We have,
\begin{align*}
\abs{d^{\Phi}(p,\ga_1(\alpha\,c_+))-d^{\Phi}(p,\ga_2(\alpha\,c_+))}
&=\abs{\int_0^{\alpha\,c_+}\Phi({\mathsf g}^t\ga_1)-\Phi({\mathsf
g}^t\ga_2)dt}
\\ &\leq K\int_0^{\alpha\,c_+} \op{dist}({\mathsf g}^t\ga_1,
{\mathsf g}^t \ga_2)^{\beta}dt.
\end{align*}
Now apply Corollary \ref{cor:ac_estimate}.
\end{proof}
\begin{lemma}  Assume $\ga_1$ and $\ga_2$ are two geodesics with
$\ga_1(\infty)=\ga_2(\infty)=\zeta$ for some $\zeta \in \pa H$.
Assume that $\rho_{\zeta}(\ga_1(0),\ga_2(0)) =0$ (i.e, $\ga_1(0)$
and $\ga_2(0)$ lie on the same horosphere around $\zeta$). Then
$$\dist({\mathsf g}^s\ga_1,{\mathsf g}^s\ga_2)\leq
\begin{cases}
 1 + d - 2\,s + \frac{s}{d} & s\leq \frac{d}{2}
\\
 e^{\left( \frac{d}{2} - s \right) }\,\left(\frac{3}{2} +\frac{s}{2}- \frac{d}{4}  \right)
  & s>
         \frac{d}{2}
\end{cases}.$$
\end{lemma}

\begin{proof} By Corollary \ref{cor:horodist}, we have
$$d(\ga_1(t),\ga_2(t))\leq 2
\sinh^{-1}\(\sinh\(\frac{d}{2}\)e^{-t}\).$$ Moreover,
$$\int_{-\infty}^\infty d({\mathsf g}^s\ga_1(t),{\mathsf g}^s\ga_2(t))e^{-\abs{t}}dt=
\int_{0}^\infty d(\ga_1(t),\ga_2(t))e^{-\abs{t-s}}dt+
e^{-s}\int_{-\infty}^0 d(\ga_1(t),\ga_2(t))e^{t}dt.$$

As before, the second integral is less than
$$  \left( 2 + c_- + d \right) \,e^{-c_- - s}\leq (2+d)e^{-s},$$ where
$c_-=\(\ga_1(-\infty)\cdot \ga_2(-\infty)\)_{\ga_1(0)}$ and
$d=d\(\ga_1(0),\ga_2(0)\)$.

Adding the first integral which can be evaluated in closed form, we
obtain the upper bound,
\begin{gather*}
2e^{-s}+ 4\,\op{csch}^{-1}\(e^s\,
         \op{csch}(\frac{d}{2})\) -
      2\,e^{s}\,\op{csch}(\frac{d}{2})\,
       \left( -1 + {\sqrt{1 + \frac{{\sinh (\frac{d}{2})}^2}{e^{2\,s}}}}
         \right)+\\
      \frac{2}{e^s}\,\( s+
      \log\(\frac{2+\sqrt{4-2e^{-2s} + 2e^{-2s}\cosh (d)}}{2 + {\sqrt{2 + 2\cosh (d)}}}\)
     \) \,\sinh (\frac{d}{2}) .
\end{gather*}
If we denote the above by $f(s)$, then $f''(s)$ is positive for
$s\geq 0$. Therefore we can overestimate $f$ for $0\leq s\leq
\frac{d}{2}$ by the line between its endpoints in this region.
Taking a power series of $g(s)=e^{s-d/2}f(s)-\(1-e^{-d}\)s$ at
$s=\infty$ we find that $g(s)$ decreases monotonically to the
constant
$$  \frac{2 + \left( 3 + \log (16) -
        2\,\log (2 + \sqrt{2 + 2\cosh (d)}) \right) \,
\sinh (\frac{d}{2})}{e^{\frac{d}{2}}}.$$
In particular, $g(s)$, and hence $e^{s-\frac{d}{2}} f(s)-s$, attains
its maximum at $s=\frac{d}{2}$. So we can overestimate $f(s)$ for
$s>\frac{d}{2}$ by
$e^{\frac{d}{2}-s}\(s-\frac{d}{2}+f\(\frac{d}{2}\)\)$ in this
region. On the other hand $f\(\frac{d}{2}\)$ is monotonically
increasing in $d$, so we can overestimate it by its limit as $d\to
\infty$ which is $ 4 - 2\,{\sqrt{5}} + 4\,\sinh^{-1}\(\frac12\) +
\sinh^{-1}\(2\)=4 - 2\,{\sqrt{5}} + \log \(\frac{29 +
13\,{\sqrt{5}}}{2}\).$

Putting these estiamtes together have:
$$2 \dist({\mathsf g}^s\ga_1,{\mathsf g}^s\ga_2)\leq \begin{cases}
  2\,\left( 1 + d - 2\,s + \frac{s\,
        \left( 2 - 2\,{\sqrt{5}} + \log \(\frac{29 + 13\,{\sqrt{5}}}{2}\)
          \right) }{d} \right) & s\leq \frac{d}{2} \\
            e^{\frac{d}{2} - s}\, \left(\,s-\frac{d}{2} + 4 -
       2\,{\sqrt{5}} + \log \(\frac{29 + 13\,{\sqrt{5}}}{2}\)
       \right)  & s> \frac{d}{2}
\end{cases}.$$
This can be further overestimated by
$$\dist({\mathsf g}^s\ga_1,{\mathsf g}^s\ga_2)\leq
\begin{cases}
 1 + d - 2\,s + \frac{s}{d} & s\leq \frac{d}{2}
\\
 e^{\left( \frac{d}{2} - s \right) }\,\left(\frac{3}{2} +\frac{s}{2}- \frac{d}{4}  \right)
  & s>
         \frac{d}{2}
\end{cases}.$$


%
\end{proof}

\medskip
\begin{corollary} Let $x, y \in H$ and $\zeta \in \pa H$. Assume that
$\rho_{\zeta}(x,y) =0$ (i.e, $x,y$ lie on the same horoball around
$\zeta$) and set $d=d\(x,y\)$.  Then
$$\int_{T}^{\infty} \op{dist}({\mathsf g}^s\ga_{x, \zeta},
{\mathsf g}^s\ga_{y, \zeta})^{\beta}ds=\begin{cases}
 \(\frac{d}{2} - T \) \,
     \(\frac{\beta}{2}+\frac{\beta\, d}{2} +1  \) +  \frac{2}{\beta}
       & T\leq
     \frac{d}{2} \\
   e^{\beta\,\left( \frac{d}{2} - T \right) }\,\(
     \frac{2}{\beta} + \frac{T}{2}-\frac{d}{4}\) & T>\frac{d}{2}
\end{cases}.$$
\end{corollary}
\medskip
\begin{proof}
Set $\ga_1=\ga_{x, \zeta}$ and $\ga_2=\ga_{y, \zeta}$. Using
$(a+x)^\beta\leq a+\beta x$ for $a\geq 1$, we estimate,

$$\dist({\mathsf g}^s\ga_1,{\mathsf g}^s\ga_2)^\beta\leq
\begin{cases}
 1 + \beta\(d - 2\,s + \frac{s}{d} \) & s\leq \frac{d}{2}
\\
 e^{\beta\,\left( \frac{d}{2} - s \right) }\,\left( \frac32
 +\frac{\beta s}2- \frac{\beta d}4  \right)    & s> \frac{d}{2}
\end{cases}.$$

This yields
$$\int_T^{\frac{d}{2}} \op{dist}({\mathsf g}^s\ga_1,{\mathsf g}^s\ga_2)^\beta=  \frac{\left( d - 2\,T \right) \,
     \left( 2\,\beta\,d^2 + 2\,\beta\,T + d\,\left( 4 + \beta - 4\,\beta\,T \right)  \right) }
     {8\,d},$$
and for $T>\frac{d}{2}$,
$$\int_T^\infty \op{dist}({\mathsf g}^s\ga_1,{\mathsf g}^s\ga_2)^\beta=
   e^{\beta\,\left( \frac{d}{2} - T \right) }\,\frac{
     4 + \beta\,\left( T- \frac{d}{2} \right)   }{2\,\beta}.$$

After overestimating the first expression, we obtain
$$\int_T^\infty \op{dist}({\mathsf g}^s\ga_1,{\mathsf g}^s\ga_2)^\beta=\begin{cases}
  \(\frac{d}{2} - T \) \,
     \(\frac{\beta}{2}\(1+  d\) +1  \) +  \frac{2}{\beta}
       & T\leq
     \frac{d}{2} \\
   e^{\beta\,\left( \frac{d}{2} - T \right) }\,\frac{
     4 + \beta\,\left( T- \frac{d}{2} \right)   }{2\,\beta}& T>
     \frac{d}{2}
\end{cases}.$$
\end{proof}

\begin{lemma}\label{lem:rho_phi_estimate}
 Let $x, y \in H$ and $\zeta \in \pa H$. Assume
that $\rho_{\zeta}(x,y) =0$ (i.e, $x,y$ lie on the same horoball
around $\zeta$).  Assume $\Phi$ is a H\"{o}lder function, with
H\"{o}lder constant $K$.  Then for all $s \geq 0$,
$$\abs{\rho_{\zeta}^{\Phi}(\ga_{x,\zeta}(s),\ga_{y,\zeta}(s))} \leq
K\,\begin{cases}
 \(\frac{d}{2} - s \) \,
     \(\frac{\beta}{2}+\frac{\beta\, d}{2} +1  \) +  \frac{2}{\beta}
       & s\leq
     \frac{d}{2} \\
   e^{\beta\,\left( \frac{d}{2} - s \right) }\,\(
     \frac{2}{\beta} + \frac{s}{2}-\frac{d}{4}\) & s>\frac{d}{2}
\end{cases}.$$
\end{lemma}

\medskip
\begin{proof} Since $\rho_{\zeta}(x,y)=0$ we have

\begin{align*}
\abs{\rho_{\zeta}^{\Phi}(\ga_{x,\zeta}(s),\ga_{y,\zeta}(s))}&=\lim_{T
\to \infty}\abs{d^{\Phi}(\ga_{y,\zeta}(s),\ga_{y,\zeta}(T))-
d^{\Phi}(\ga_{x,\zeta}(s),\ga_{x,\zeta}(T))}=
\\ & = \lim_{T \to \infty}\int_{s}^{T} \abs{\Phi({\mathsf g}^t\ga_{y,\zeta}) -
\Phi({\mathsf g}^t\ga_{x,\zeta})}dt \\ &\leq K\int_{s}^{\infty}
\,\op{dist}({\mathsf g}^t \ga_{y,\zeta}, {\mathsf
g}^t\ga_{x,\zeta})^{\beta}\,dt.
\end{align*}
\end{proof}

 For any two points $x \in H$ and $\xi \in  \pa H$, let
$\ga_{x,\xi}(t)$ be a geodesic connecting $x$ and $\xi$, i.e.,
$\ga_{x, \xi}(0)=x$ and $\lim_{t \to \infty} \ga_{x, \xi}(t) = \xi$.

Now recall that any CAT($-\kappa$) space is $2\delta$-hyperbolic for
some $\delta>0$ depending only on $\kappa$. For the remainder of
this section, let $\delta$ be the smallest such constant for $H$.
\begin{lemma}\label{lem:hyperbolic_laws}
Suppose $p,q\in H$ and $\zeta,\nu\in\pa H$ are chosen so that
$(\zeta\cdot\nu)_p\geq d(p,q)$ and $\rho_\zeta(p,q)\geq 0$, then
$$d(\ga_{p,\zeta}( (\zeta\cdot\nu)_p ),\ga_{p,\nu}( (\zeta\cdot\nu)_p ))\leq 2\delta,$$
$$d(\ga_{q,\zeta}( (\zeta\cdot\nu)_p +
\rho_{\zeta}(p,q)),
 \ga_{q,\nu}( (\zeta\cdot\nu)_p +\rho_{\nu}(p,q))) \leq 4\delta.$$
$$\max\{d(p,\ga_{q,\zeta}(\rho_{\zeta}(p,q))),d(p,\ga_{q,\nu}(\rho_{\nu}(p,q)))\}
\leq d(p,q)+2\delta$$
\end{lemma}

\begin{proof}  The first inequality
follows directly from writing down the conditions for
$\delta$-hyperbolicity. For $(\zeta\cdot\nu)_p\geq d(p,q)$ we have
$\abs{\rho_{\zeta}(p,q) -\rho_{\nu}(p,q)}\leq 2\delta.$

Since $\rho_{\zeta}(p,q)\geq 0$, we therefore have
$\rho_{\nu}(p,q)\geq -2\delta$. By Lemma \ref{lem:increasing_dist},
$d(p,\ga_{q,\zeta}(\rho_{\zeta}(p,q))) \leq d(p,q)$ and hence
$d(p,\ga_{q,\nu}(\rho_{\nu}(p,q))) \leq d(p,q)+2\delta$ giving the
second inequality.

Finally, observe that \begin{align*} \max\{ (\zeta\cdot\nu)_p
+\rho_{\zeta}(p,q),
 (\zeta\cdot\nu)_p +\rho_{\nu}(p,q)\} \leq
(\zeta\cdot\nu)_p+\frac{\rho_{\zeta}(p,q)+\rho_{\nu}(p,q)}{2}+\delta
=(\zeta\cdot\nu)_p+\delta
\end{align*}
Thus we have
\begin{align*}
d(\ga_{q,\zeta}( (\zeta\cdot\nu)_p + \rho_{\zeta}(p,q)),
 \ga_{q,\nu}( (\zeta\cdot\nu)_p +\rho_{\nu}(p,q))) \leq
d(\ga_{q,\zeta}((\zeta\cdot\nu)_q),\ga_{q,\nu}((\zeta\cdot\nu)_q)+2\delta)
\leq 4\delta.\end{align*}
\end{proof}

\medskip
At this point we need to make an assumption about the way we extend
geodesics backwards.
\begin{definition}\label{def:flip}
We say that $H$ has {\em bounded flip} if there exist positive
constants $R_1$ and $R_2$ such that for all $\zeta,\nu\in\pa H$ and
$p\in H$, it is possible to choose two geodesics $\ga_1$ and $\ga_2$
such that $\ga_1^+=\zeta$, $\ga_2^+=\nu$, $\ga_1(0)=p=\ga_2(0)$ and
we have
$$(\ga_1^-\cdot \ga_2^-)_{\ga_1(0)} \geq R_1(\ga_1^+\cdot \ga_2^+)_{\ga_1(0)}-R_2.$$
\end{definition}

\begin{remark}
Observe that the Cayley graphs of free groups have this property
(with $R_1>0$ arbitrary) since the backward endpoints of geodesics
through $p$ can be chosen independently of the forward points. For a
simply connected manifold $M$ with pinched negative curvature
between $-1$ and $-b^2$, we have $b (\ga_1^-\cdot
\ga_2^-)_{\ga_1(0)}\geq (\ga_1^+\cdot \ga_2^+)_{\ga_1(0)}$. In
particular, $M$ has bounded flip. On the other hand, it is easy to
construct simply connected manifolds with unbounded negative
curvature less than  $-1$ which do not possess bounded flip. Simply
take a geodesic $\ga$ through $p$ whose forward ray passes through a
region of curvature $-1$, but whose backward ray has the property
that every 2-plane tangent to $\ga'(-t)$ has sectional curvature
$-e^{2t}$. Therefore for any sequence of geodesics $\ga_i$ through
$p$ with $\ga_i\to \ga$ we have $\frac{(\ga_i^+\cdot
\ga^+)_{p}}{(\ga_i^-\cdot \ga^-)_{p}}\to \infty$.
\end{remark}

The next proposition shows that CAT(-1) spaces which are either tree
like or have pinched curvatures in a weak sense and admit a
cocompact group of isometries have bounded flip. This family
includes other important examples such as hyperbolic buildings.
Conjecturally this any CAT(-1) space with a cocompact group of
isometries has bounded flip. However, we were not able to show this.

\begin{proposition}
Suppose $G$ acts cocompactly on $H$, and there is a $K<-1$ such that
every triangle with two endpoints in $\pa H$ either
\begin{enumerate}
\item has both sides meeting the interior vertex along a common segment, or
\item  is fatter than its comparison triangle in $\HH^2_K$ the
plane of  constant curvature $K$.
\end{enumerate}
Then $H$ has bounded flip.
\end{proposition}
\begin{proof}
Since the assumptions imply that $H$ is a proper metric space, even
if it is not geodesically complete, any Cauchy sequence of geodesic
segments with respect to Hausdorff distance is convergent. This
implies that $\pa H$ is complete, and hence compact, with respect to
the quasimetric $\pi_p(\zeta,\nu)=e^{-(\zeta,\nu)_p}$. In
particular, a sequence of geodesic rays through a fixed point has a
subsequence which converges to a ray.

Therefore if $H$ does not have bounded flip, then we may assume
there is a fixed geodesic $\ga$ and a sequence of geodesics $\ga_i$
with $\ga_i(0)=\ga(0)$ and $\ga_i\to\ga$ such that
$\lim_{i\to\infty}\frac{(\ga_i^+\cdot \ga^+)_{\ga(0)}}{(\ga_i^-\cdot
\ga^-)_{\ga(0)}}\to \infty$ where $(\sigma^-\cdot
\ga^-)_{\ga(0)}\leq (\ga_i^-\cdot \ga^-)_{\ga(0)}$ for any other
geodesic $\sigma$ satisfying $\sigma(0)=\ga(0)$ and
$\sigma^+=\ga_i^+$.

For triangles of the first type the geodesics from $\ga(0)$ to
$\ga_i^+$ and $\ga^+$ can be extended in the opposite direction by a
common geodesic, so that there is no constraint imposed on the
choice of constant $R_1$.

For a triangle in the hyperbolic plane of constant curvature
$K=-k^2$ with one point at $p\in \HH^2$ and the other two vertices
at $\xi,\zeta\in \pa \HH^2$, the angle $\theta$ at $p$ is given by
$\theta=\arcsin\(2 e^{-k(\xi,\zeta)_p}\)$. Therefore, the second
assumption implies that $(\ga_i^-,\sigma^-)_{\sigma(0)}\geq
\frac{1}{k}(\ga_i^+,\sigma^+)_{\sigma(0)}$. So that $R_1$ may be
chosen to be $\frac{1}{k}$.
\end{proof}

%

\begin{proposition}\label{prop:RN_holder} Assume $H$ has bounded flip. For $p \in H$ we
define $\pi_p(\xi,\zeta)=e^{-(\xi\cdot\zeta)_p}$.
 There exists $\ep>0$ and constant $D_\Phi$  depending only on
$\Phi$ such that for any $q \in H$ we have
$$ \sup_{\pi_p(\zeta, \nu) \leq e^{- d(p,q)}}
\frac{\abs{\rho_{\zeta}^{\Phi}(p,q)
-\rho_{\nu}^{\Phi}(p,q)}}{\pi_p(\zeta,\nu)^{\ep}} <D_{\Phi}e^{\ep
d(p,q)}.$$
\end{proposition}

\begin{proof}

Let $f(p,q)$ be the left hand side of the expression. Even though,
$\rho_\zeta^{\Phi}(p,q) = -\rho_\zeta^{\Phi}(q,p)$ for all
$\zeta\in\partial H$, we may still have $f(p,q)\neq f(q,p)$ since
$(\zeta\cdot \nu)_p=(\zeta\cdot
\nu)_q+\frac12\(\rho_\zeta(q,p)+\rho_\nu(q,p)\)$. However, by
cocompactness of $\Gamma$, there is a $\ga\in\Gamma$ such that
$d(\ga p,q)<\op{diam}(H/\Ga)$ and $\pi_p(\zeta,\nu)=\pi_{\ga p}(\ga
\zeta,\ga \nu)$. The lower semi-continuity of $f(p,q)$ implies that
$f(p,q)\leq f(q,p)e^{2 \eps \op{diam}(H/\Ga)}$. Therefore, without
loss of generality we assume that $\rho_{\zeta}(p,q)\geq 0$.

Set $r = \pi_p(\zeta, \nu)\leq e^{- d(p,q)}$ and $c_+
=(\zeta\cdot\nu)_p= -\log{(r)} \geq d(p,q)$ (so $e^{-c_+} = r$).

From the definition of $\rho^\Phi$, for all $x, y \in H$ and $\xi
\in \pa H$ and $t_1, t_2 \geq 0$, we have
$$\rho_{\xi}^{\Phi}(x,y) = \rho_{\xi}^{\Phi}(\ga_{x, \xi}(t_1),
\ga_{y,\xi}(t_2)) + d^{\Phi}(y,\ga_{y,\xi}(t_2)) -
d^{\Phi}(x,\ga_{x,\xi}(t_1)).$$ And points $p$ and
$\ga_{q,\xi}(\rho_{\xi}(p,q))$ lie on the same horoball for all
$p,q\in H$ $\xi \in \pa H$.

 Now let us estimate $\abs{\rho_{\zeta}^{\Phi}(x,y)
-\rho_{\nu}^{\Phi}(x,y)}$. Let $a\in (0,1)$ be a parameter which we
will specify later.

We will concentrate on $4$ points
$$\ga_{p,\zeta}(\alpha\, c_+ ),  \ga_{q,\zeta}(\alpha\, c_+ +
\rho_{\zeta}(p,q)), \text{ and } \ga_{p,\nu}(\alpha\, c_+ ),
 \ga_{q,\nu}(\alpha\, c_+ +\rho_{\nu}(p,q)).$$
 It is not difficult to see that the first two points and the last two points lie on
the same horoball.

Now we can make the following estimate.
\begin{align*}
\rho_{\zeta}^{\Phi}(p,q) -\rho_{\nu}^{\Phi}(p,q) = &
\rho_{\zeta}^{\Phi}(\ga_{p,\zeta}(\alpha\, c_+ ),
\ga_{q,\zeta}(\alpha c_+ +\rho_{\zeta}(p,q))) -
\rho_{\nu}^{\Phi}(\ga_{p,\nu}(\alpha c_+), \ga_{q,\nu}(\alpha\,c_+
+\rho_{\nu}(p,q)))\\&+ d^{\Phi}(q,\ga_{q,\zeta}(\alpha\,c_+
+\rho_{\zeta}(p,q))) - d^{\Phi}(q,\ga_{q,\nu}(\alpha\,c_+
+\rho_{\nu}(p,q)))\\& -d^{\Phi}(p,\ga_{p,\zeta}(\alpha\,c_+
))+d^{\Phi}(p,\ga_{p,\nu}(aT))
\end{align*}

We will estimate each line separately. However, we will do the
estimates in the order they appear.

The estimates of each $\rho^{\Phi}$ term are similar and we will
consider them together. While Lemma \ref{lem:hyperbolic_laws} tells
us that $d(p,\ga_{q,\zeta}(\rho_{\zeta}(p,q)))\leq d(p,q)+2\delta$,
we need a better bound when $d(p,q)$ is small. Using Corollary
\ref{cor:horodist}, we can solve to also show that
$d(p,\ga_{q,\zeta}(\rho_{\zeta}(p,q)))\leq
\cosh^{-1}\(\frac12\(\cosh(d(p,q))^2+1\)\).$  Since
$f(t)=\frac{\cosh^{-1}\(\frac12\(\cosh(t)^2+1\)\)}{t}$ is monotone
and $\lim_{t\to 0}f(t)=1$ and $\lim_{t\to \infty}f(t)=2$, we have
$\frac{d(p,\ga_{q,\zeta}(\rho_{\zeta}(p,q)))}{d(p,q)}\leq
\min\set{f(d(p,q)),1+\frac{2\delta}{d(p,q)}}\leq C_0$, where
$C_0\in(1,2)$ is a constant depending only on $\delta$.

Since $c_+>d(p,q)$, if we choose any $1>a >\frac{C_0}{2}$ then $a
c_+>\frac{d(p,\ga_{q,\zeta}(\rho_{\zeta}(p,q)))}{2}$ and by
symmetry, $a c_+>\frac{d(p,\ga_{q,\nu}(\rho_{\nu}(p,q)))}{2}$.
Therefore by Lemma \ref{lem:rho_phi_estimate} we have,
\begin{align*}
&\max\set{\abs{\rho_{\zeta}^{\Phi}(\ga_{p,\zeta}(\alpha\, c_+ ),
\ga_{q,\zeta}(\alpha\, c_+ +\rho_{\zeta}(p,q)))},
\abs{\rho_{\nu}^{\Phi}(\ga_{p,\nu}(\alpha\, c_+
),\ga_{q,\nu}(\alpha\, c_+ +\rho_{\nu}(p,q)))}}
\\ & \leq 2 Ke^{-\beta \alpha\, c_+ }\max\left\{ e^{\beta\left(
\frac{d(p,\ga_{q,\zeta}(\rho_{\zeta}(p,q)))}{2}\right)},e^{\beta
\left(\frac{d(p,\ga_{q,\zeta}(\rho_{\nu}(p,q)))}{2}\right)}\right\}\(\frac{a
c_+}{2} +\frac{2}{\beta}\)
\\ & \leq
K e^{\beta \alpha\, c_+ }e^{\beta
\left(\frac{d(p,q)}{2}+\delta\right)}\(\frac{\alpha\, c_+}{2}
+\frac{2}{\beta}\) \leq \frac{2 K e^{\beta
\left(\frac{d(p,q)}{2}+\delta\right)}}{\beta} e^{-\frac{3}{4}\beta a
c_+ } \\
&\leq \frac{2 K e^{\beta \delta}}{\beta}e^{\frac{\beta}{2} c_+}
r^{\frac{3}{4}\beta \alpha\,}=\frac{2 K e^{\beta \delta}}{\beta}
r^{\beta\(\frac{3}{4}\alpha\,-\frac{1}{2}\)}.
\end{align*}

Here we have used the fact that $e^{-x}(A+Bx) \leq A
e^{-(1-\frac{B}{A})x}$, for $x \geq 0$.

Observe that $\ga_{p,\zeta}(0)=\ga_{p,\nu}(0)=p$, so we can apply
Lemma \ref{lem:dist_difference}
\begin{align*}
&\abs{d^{\Phi}(p,\ga_{p,\zeta}(\alpha\, c_+
))-d^{\Phi}(p,\ga_{p,\nu}(\alpha\, c_+
))}\\
&\leq K\,\(\(\frac{2}{\beta} + \frac{c_-}{2}\)\,e^{-\beta\,c_-}\,
       +\(\frac{2}{\beta}+\frac{(1-\alpha)c_+}{2}\)e^{-\beta \,(1-\alpha)\,c_+}\)\\
& \leq \frac{2K}{\beta}\left(e^{-\frac{3}{4}\beta
c_-}+e^{-\frac{3}{4}\beta (1-\alpha\,)c_+}\right)\\
\intertext{recalling that $c_- \geq R_1 c_+-R_2$ yields,}
 & \leq \frac{2K}{\beta}\left(e^{-\frac{3}{4}\beta
\(R_1 c_+ -R_2\)}+e^{-\frac{3}{4}\beta (1-\alpha\,)c_+}\right) \\
& =\frac{2K}{\beta}\left(e^{\frac{3}{4}\beta R_2}r^{\frac{3}{4}\beta
R_1}+r^{\frac{3}{4}\beta (1-\alpha\,)}\right) .
\end{align*}

Now we will repeat this to the second estimate at the point $q$.
\begin{align*}
&\abs{d^{\Phi}(q,\ga_{q,\zeta}(\alpha\, c_+ +\rho_{\zeta}(p,q))) -
d^{\Phi}(q,\ga_{q,\nu}(\alpha\, c_+ +\rho_{\nu}(p,q)))}  \\ & \leq
\abs{d^{\Phi}(q,\ga_{q,\zeta}(\alpha\, c_+ +\rho_{\zeta}(p,q))) -
d^{\Phi}(q,\ga_{q,\nu}(\alpha\, c_+ +\rho_{\zeta}(p,q)))} \\ &+
\abs{d^{\Phi}(\ga_{q,\nu}(\alpha\, c_+
+\rho_{\zeta}(p,q)),\ga_{q,\nu}(\alpha\, c_+ +\rho_{\nu}(p,q)))}
\intertext{Note that $(\zeta\cdot
\nu)_q=c_++\frac12\(\rho_\zeta(p,q)+\rho_\nu(p,q)\)$, so $\abs{
c_++\rho_\zeta(p,q)-(\zeta\cdot\nu)_q}\leq2\delta$. We use Lemma
\ref{lem:dist_difference} (with $(\zeta\cdot\nu)_q$ in place of
$c_+$) again to estimate the first term, and estimating the second
term explicitly, we have } \leq &
\frac{2K}{\beta}\left(e^{-\frac{3}{4}\beta
(\ga_{q,\zeta}(-\infty)\cdot\ga_{q,\nu}(-\infty))_q}+e^{-\frac{3}{4}\beta
\(1-\alpha\,(1+\frac{2\delta}{(\zeta\cdot\nu)_q})-
(1-\alpha\,)\frac{\rho_\zeta(p,q)}{(\zeta\cdot\nu)_q}\)(\zeta\cdot\nu)_q}\right)\\
&+K\abs{\rho_{\zeta}(p,q)-\rho_{\nu}(p,q)}\\
\leq & \frac{2K}{\beta}\left(e^{-\frac{3}{4}\beta (R_1c_+-R_2)}+
e^{\frac32\beta\delta}e^{-\frac{3}{4}\beta
(1-\alpha\,)\((\zeta\cdot\nu)_q-\rho_\zeta(p,q)\)}\right)\\
&+K\abs{\rho_{\zeta}(p,q)-\rho_{\nu}(p,q)}\\
 \leq &\frac{2K}{\beta} \left(e^{\frac34\beta R_2}r^{\frac34\beta
R_1} +e^{3\beta \delta}r^{\frac34\beta(1-\alpha\,)}\right)+
K\abs{\rho_{\zeta}(p,q)-\rho_{\nu}(p,q)}.
\end{align*}

Now let us specify $\ep$ and $a$. Suppose $\ep_0$ is the least upper
bound of $\ep$ for which $\pi_p^{\ep}$ is a metric. In a moment we
will need to choose $\ep$  so that
$$\ep\leq  \min{\left(\ep_0,\frac34 \beta \alpha\,,\frac34
\beta R_1,\frac34 \beta (1-\alpha\,),\frac34\beta
\alpha\,-\frac{\beta}{2} \right)}.$$ Recall
$1>a>\frac{C_0}{2}>\frac12$. So the best upper bound for $\ep$ will
be when $\frac34(1-\alpha\,)=\frac34 \alpha\,-\frac12$ or
$\alpha=\frac56$, assuming this is larger than $C_0$. This is the
case when $\delta$ is chosen for a CAT($-1$) space. This yields an
upper bound of $\ep <\min\set{\ep_0,\frac34\beta
R_1,\frac{\beta}{8}}.$
\begin{align*}
&\frac{\abs{\rho_{\zeta}^{\Phi}(p,q)-\rho_{\nu}^{\Phi}(p,q)}}{r^{\ep}}
  \\
&\leq \frac{2 K e^{\beta \delta}}{\beta}
r^{\beta\(\frac{3}{4}\alpha\,-\frac{1}{2}\)-\ep}
+\frac{2K}{\beta}\left(e^{\frac{3}{4}\beta R_2}r^{\frac{3}{4}\beta
R_1-\ep}+r^{\frac{3}{4}\beta (1-\alpha\,)-\ep}\right) \\
&+\frac{2K}{\beta} \left(e^{\frac34\beta R_2}r^{\frac34\beta
R_1-\ep} +e^{3\beta \delta}r^{\frac34\beta(1-\alpha\,)-\ep}\right)+
K\frac{\abs{\rho_{\zeta}(p,q)-\rho_{\nu}(p,q)}}{r^{\ep}}.
\end{align*}
The last term can be controlled since $\sup_{\pi_p(\zeta,\nu)\leq
e^{-d(p,q)}} \frac{\abs{\rho_{\zeta}(p,q)-\rho_{\nu}(p,q)}}{r^{\ep}}
\leq Ce^{\ep d(p,q)}$. Since all powers of $r$ are positive and $r
\leq e^{-d(p,q)}\leq 1$ we have finished the proof. \end{proof}.
\newpage
\section{Decaying condition.}\label{app:decay_cond}
Recall that $S_p(q,r)\in \pa H$ is the shadow of the ball $B(q,r)
\ssu H$ with respect to point $p\in H$.

\begin{theorem}[Shadow Lemma
\cite{Kaimanovich04}]\label{thm:shadow_lemma} Assume $G$ acts
cocompactly on $H$. If $\set{\mu_p^\Phi}$ is a Gibbs stream
associated to $\Phi$, then for any $r>0$ and $p\in H$ there exists
$C_{\Phi}(r)>0$ such that
$$\frac{1}{C_{\Phi}(r)} e^{[-d^{\Phi}(p,q)- \la_{\Phi} d(p,q)]}\leq \mu_p^\Phi(S_p(q,r))
\leq C_{\Phi}(r) e^{[-d^{\Phi}(p,q)- \la_{\Phi} d(p,q)]},$$ for all
$p,q \in H$.
\end{theorem}

\begin{remark}
The proof of the above theorem shows a little more. If $G<\Isom(H)$
is arbitrary, then the same conclusion holds but with $C_\Phi$ also
depending on $p$ and $d(q,G\cdot p)$.
\end{remark}

Recall that under our assumptions on $\Phi$, we have $\la_{\Phi}=0$.

Let $\{m_p\}_{p\in H}$ be a Hausdorff stream, i.e., any Gibbs stream
associated to $\Phi=0$. Let $\la_0$ be its critical constant. We
have the following

\begin{corollary}\label{cor:shadow}  There exists a constant $K_{\Phi}$ such
that for any $\Phi$ satisfying the normalization $\lambda_\Phi=0$,
$$\frac{1}{K_{\Phi}} e^{-\la_0s} \leq \int_{\pa H}
e^{-d^{\Phi}(p,\ga_{p,\xi}(s))}dm_p(\xi) \leq K_{\Phi}e^{-\la_0s},$$
for every $s >0$.
\end{corollary}

\medskip
\begin{proof}
Fix $r>0$.  Take the largest non-intersecting collection of balls
$\{B(q_i, r/2)\}_{i=1}^k$ such that $d(p,q_i)=s$.

{\bf Claim:} {\it The collection of $\{S_p(q_i, r/4)\}_{i=1}^k$ is
disjoint}

\medskip

$\triangle$ Assume that $\xi \in S_p(q_1, r/4) \cap S_p(q_2,r/2)$.
Let $t_i$, $i=1,2$ are such that $\ga_{p,\xi}(t_i) \in B(q_i, r/4)$,
for $i=1,2$. It is clear that $s-r/4 < t_1, t_2 < s_r/4$. Thus $
d(q_i, \ga_{p,\xi}(s)) < r/2$ for $i=1,2$.  Therefore,
$B(q_1,r/2)\cap B(q_2,r/2) \not= \emptyset$. Contradiction.
$\blacktriangle$

 An easy argument shows that $\cup_{i=1}^k B(q_i, r)$
cover a sphere of radius $s$ around $p$.

{\bf Claim:} {\it $\{S_p(q_i, r)\}_{i=1}^k $ is a Besecovitch cover
for $\pa H$ with Lebesgue number $B_{\Phi}(r)=C_{\Phi}(9r/4)
C_{\Phi}(r/4)$.}

\medskip
$\triangle$ Let $\xi \in \pa H$. Assume that $\xi \in S_p(q_i, r)$
then easy calculation shows that $d(q_i, \ga_{p,\xi}(s))\leq 2r$.
Thus $B(q_i, r/4) \ssu B(\ga_{p,\xi}(s), 9r/4 )$ and therefore,
$$S_p(q_i, r/4) \ssu S_{p}(\ga_{p,\xi}(s), 9r/4).$$
Since $\{S_p(q_i,r/4)\}$ are disjoint, by Theorem above there are at
most $C_{\Phi}(9r/4 ) C_{\Phi}(r/4)$ Shadows that cover $\xi$. This
proves the Claim. $\blacktriangle$

\medskip
Let $\{S_p(q_i, r)\}_{i=1}^k$ be a Besecovitch cover as above. By
choice of the $q_i$ we have $c_+=(\zeta\cdot\nu)_p\geq s-r$. By
using the proof of lemma \ref{lem:dist_difference} with the last
estimate in Lemma \ref{lem:integral_estimate} (with $T=s=c_+ + r$)
for any two points $\xi, \nu \in S_p(q,r)$ we have
$$ \abs{d^{\Phi}(p,
\ga_{p,\xi}(s))-d^{\Phi}(p, \ga_{p,\nu}(s))} \leq
K\(\frac{5}{\beta}+r^{1+\beta}\),$$ where $K$ and $\beta$ are the
H\"{o}lder constant and exponent of $\Phi$ respectively.

So $$\aligned \int_{\pa H} e^{-d^{\Phi}(p,\ga_{p,\xi}(s))}dm_p(\xi)&
\leq \sum_{i=1}^k
\int_{S_p(q_i,r)} e^{-d^{\Phi}(p,\ga_{p,\xi}(s))}dm_p(\xi)  \\
& \leq \sum_{i=1}^k \int_{S_p(q_i,r)} e^{-d^{\Phi}(p,q_i)}
e^{K\(\frac{5}{\beta}+r^{1+\beta}\)}dm_p(\xi)
\\& = \sum_{i=1}^k e^{-d^{\Phi}(p,q_i)}
e^{K\(\frac{5}{\beta}+r^{1+\beta}\)}m_p(S_p(q_i,r))  \\& \leq
\sum_{i=1}^k C_{\Phi}(r) \mu_p^\Phi(S_p(q_i,r))
e^{K\(\frac{5}{\beta}+r^{1+\beta}\)}m_p(S_p(q_i,r))
\\ & \leq \sum_{i=1}^k C_{\Phi}(r) \mu_p^\Phi(S_p(q_i,r))
e^{K\(\frac{5}{\beta}+r^{1+\beta}\)} C_{0}(r) e^{-\la_0s}  \\ & \leq
e^{K\(\frac{5}{\beta}+r^{1+\beta}\)} C_{0}(r) e^{-\la_0s} B_{\Phi}
C_{\Phi}(r)
\endaligned$$
 Similarly we can obtain a lower bound
$$\aligned \int_{\pa H} e^{-d^{\Phi}(p,\ga_{p,\xi}(s))}dm_p(\xi)&
\geq \frac{1}{B_{\Phi}}\sum_{i=1}^k \int_{S_p(q_i,r)}
e^{-d^{\Phi}(p,\ga_{p,\xi}(s))}dm_p(\xi)  \\ & \geq
\frac{1}{B_{\Phi}}\sum_{i=1}^k \int_{S_p(q_i,r)}
e^{-d^{\Phi}(p,q_i)} e^{-K\(\frac{5}{\beta}+r^{1+\beta}\)}dm_p(\xi)
\\& = \frac{1}{B_{\Phi}}\sum_{i=1}^k e^{-d^{\Phi}(p,q_i)}
e^{-K\(\frac{5}{\beta}+r^{1+\beta}\)}m_p(S_p(q_i,r)) \\& \geq
\frac{1}{B_{\Phi}} \sum_{i=1}^k \frac{1}{C_{\Phi}(r)}
\mu_p^\Phi(S_p(q_i,r))
e^{-K\(\frac{5}{\beta}+r^{1+\beta}\)}m_p(S_p(q_i,r))
\\ & \geq \frac{1}{B_{\Phi}} \sum_{i=1}^k \frac{1}{C_{\Phi}(r)}
\mu_p^\Phi(S_p(q_i,r)) e^{-K\(\frac{5}{\beta}+r^{1+\beta}\)}
\frac{1}{C_{0}(r)} e^{-\la_0s}  \\
& \geq \frac{1}{B_{\Phi}} e^{-K\(\frac{5}{\beta}+r^{1+\beta}\)}
\frac{1}{C_{0}(r)C_{\Phi}(r) } e^{-\la_0s}
\endaligned$$
Now since the integral does not depend on $r$ we can find $K_{\Phi}$
satisfying the lemma. This proves the Lemma. \end{proof}

\begin{remark}
While we will not need this, the above proof generalizes to show
that there exists a constant $K=K_{\Phi_1,\Phi_2}$ such that for any
H\"{o}lder functions $\Phi_1,\Phi_2$, not necesarily normalized,
$$\frac{1}{K}  \leq e^{(\la_{\Phi_1}+\la_{\Phi_2})s}\int_{\pa H}
e^{-d^{\Phi_1}(p,\ga_{p,\xi}(s))}d\mu^{\Phi_2}_p(\xi) \leq K,$$ for
every $s >0$.
\end{remark}

\begin{proposition}\label{prop:positive_average}
For any fixed $p\in H$, the function $\Phi+\la_\Phi$ has a positive
geodesic average along any geodesic which stays within a fixed
distance from the orbit $G\cdot p$. If $G$ acts cocompactly on $H$,
then for some $\ep>0$,  $\Phi+\la_\Phi$ has geodesic average $\ep$
at all $p\in H$.
\end{proposition}

\begin{proof}
First observe that for any $p,q \in H$ and $U \in \pa H$ we have
that
\begin{align*}
\mu_q^\Phi(U) =& \int_U d\mu_q^\Phi(\xi)  = \int_U \frac{d
\mu_q^\Phi}{d \mu_p^\Phi }(\xi) d \mu_p^\Phi(\xi) \leq \int_U
e^{(\|\Phi\|_{\infty} + \la_{\Phi})d(p,q)} d \mu_p^\Phi(\xi) = \\& =
e^{(\|\Phi\|_{\infty} + \la_{\Phi})d(p,q)} \mu_p^\Phi(U).
\end{align*}

By Theorem \ref{thm:shadow_lemma} we have that for any fixed $R>0$,
$$\frac{1}{C_{\Phi}(r)} e^{[-d^{\Phi}(p,q)- \la_{\Phi} d(p,q)]}\leq
\mu_p^\Phi(S_p(q,r)) \leq C_{\Phi}(r) e^{[-d^{\Phi}(p,q)- \la_{\Phi}
d(p,q)]},$$ for all $q$ with $d(q,G\cdot p)<R$.

Fix $r$ and set the constant $C_0=C_{\Phi}(r)$. In particular, for
any $K>1$ and corresponding $\ep=\frac{1}{K C_0}>0$, there exists
$T_p$ such that for all $q$ with the property $d(p,q)\geq T_p$ and
$d(q,G\cdot p)<R$, we have $\mu_p^\Phi(S_p(q,r)) \leq \ep$. This
implies that
$$d^{\Phi}(p,q)+ \la_{\Phi} d(p,q) \geq -\log(\ep) - \log(C_0) = \log(K)>0.$$
This proves the first statement. Since it is not clear that $\ep$
varies continuously in $p$, we will have to do a little extra work
for the stronger statement when $G$ acts cocompactly.

For all $p^{\pr}$ such that $d(p, p^{\pr}) \leq \op{diam}(H/G)$, by
corollary \ref{cor:compare_geodesics}
$$|d^{\Phi}(p,q) - d^{\Phi}(p^{\pr},q)| \leq D(\op{diam}(H/G)).$$

Thus for all $q$ such that $d(p^{\pr}, p) \leq \op{diam}(H/G)$ and
$d(p, q) \geq T_p$ we have
\begin{align*}
d^{\Phi}(p^{\pr},q)+ \la_{\Phi} d(p^{\pr},q) &\geq d^{\Phi}(p,q)+
\la_{\Phi} d(p,q) - D(\op{diam}(H/G)) - \la_{\Phi} d(p, p^{\pr})
\\ & \geq \log(K) -D(\op{diam}(H/G)) - \la_{\Phi} d(p,
p^{\pr}).\end{align*}

Take $K$ such that $a=\log(K) -D(\op{diam}(H/G)) - \la_{\Phi} d(p,
p^{\pr})$ is positive.

By $G$ invariance we have proved that for all $p,q$ such that
$d(p,q) \geq T=T_p + \op{diam}(H/G)$ we have $$d^{\Phi}(p,q)+
\la_{\Phi} d(p,q)\geq a.$$

Now for every geodesic $\ga$ and $N\geq T$ we have

\begin{align*}
&\int_{0}^N \Phi(\dot{\ga}(t)) dt + \la_{\Phi} N = \\ &=
\left(\sum_{i=0}^{[N/T]-2} \int_{iT}^{(i+1)T} \Phi(\dot{\ga}(t)) dt
+ \la_{\Phi} T \right) + \int_{[N/T]T-T}^{N} \Phi(\dot{\ga}(t)) dt +
\la_{\Phi} (N-[N/T]T-T)  \\ & = \left(\sum_{i=0}^{[N/T]-2}
d^{\Phi}(\ga(it), \ga((i+1)T)+ \la_{\Phi}d(\ga(it), \ga((i+1)T)
\right)  \\ &\phantom{\hspace{2cm}}+ d^{\Phi}(\ga([N/T]T-T),
\ga(N))+ \la_{\Phi}d(\ga([N/T]T-T), \ga(N))\\ & \geq ([N/T]-1)a + a
= [N/T] a \geq N \frac{a}{T} - a
\end{align*}

This shows uniformly positive average along all geodesics.
\end{proof}


\end{document}